\documentclass[a4paper,twoside]{article}
\usepackage{amssymb}
\usepackage{amsbsy}
\usepackage{amsmath}
\usepackage{theorem}
\usepackage{epsfig}
\usepackage{graphicx}

\evensidemargin\oddsidemargin 
\numberwithin{equation}{section}
\newtheorem{theorem}{Theorem}[section]
\newtheorem{corollary}[theorem]{Corollary}
\newtheorem{conjecture}[theorem]{Conjecture}
\newtheorem{lemma}[theorem]{Lemma}
\newtheorem{proposition}[theorem]{Proposition}
{\theorembodyfont{\rmfamily}
\newtheorem{definition}[theorem]{Definition}
}
{\theorembodyfont{\rmfamily}
\newtheorem{remark}[theorem]{Remark}
}
{\theorembodyfont{\rmfamily}

}
{\theorembodyfont{\rmfamily}

}
{\theorembodyfont{\rmfamily}
\newtheorem{examples}[theorem]{Examples}
}
\begin{document}

\title{\textbf{On the String-Theoretic Euler Number }\\
\textbf{of a Class of Absolutely Isolated Singularities}}
\author{\thanks{{\scriptsize GSRT-fellow supported by the E. U. and the Greek
Ministry of Research and Technology.}\medskip \newline
\medskip\ \newline
}\ \textbf{Dimitrios I. Dais} \\
%EndAName
\noindent{\scriptsize Mathematics Department, Section of Algebra and
Geometry, University of Ioannina}\\
{\scriptsize GR-45110 Ioannina, Greece, E-mail: ddais@cc.uoi.gr}}
\date{}
\maketitle

\begin{abstract}
\noindent An explicit computation of the so-called string-theoretic
E-function{\small \ of a normal complex variety $X$ with at most
log-terminal singularities can be achieved by constructing one
snc-desingularization of $X$, accompanied with the intersection graph of the
exceptional prime divisors, and with the precise knowledge of their
structure. In the present paper, it is shown that this is feasible for the
case in which $X$ is the underlying space of a class of absolutely isolated
singularities (including both usual $\mathbf{A}_{n}$-singularities and
Fermat singularities of arbitrary dimension). As byproduct of the exact
evaluation of $\ e_{\text{str}}\left( X\right) =$ lim$_{u,v\rightarrow
1}\,E_{\text{str}}\left( X;u,v\right) $, for this class of singularities,
one gets counterexamples to a conjecture of Batyrev concerning the
boundedness of the string-theoretic index. Finally, the string-theoretic
Euler number is also computed for global complete intersections in $\mathbb{P%
}_{\mathbb{C}}^{N}$ with prescribed singularities of the above type. }
\end{abstract}

%====ABSTRACT===============================

%=============FIRST SECTION=======================

\section{Introduction\label{INTRO}}

\noindent {}The so-called $E_{\text{str}}$-\textit{polynomials }$E_{\text{str%
}}\left( X;u,v\right) $ of normal complex varieties $X$ with at most
Gorenstein quotient or toroidal singularities were introduced in \cite{BD},
and were used as main tools in \cite{BD} and \cite{BB} for the proof of
several mirror-symmetry identities. More recently, Batyrev \cite{BATYREV1}
generalized this notion also for $X$'s having at most log-terminal
singularities, by introducing $E_{\text{str}}$-\textit{functions} instead
which may be not even rational. These new invariants have already found lots
of applications in the study of log-flips and of cohomological McKay
correspondence. (See \cite[1.6, 4.11 and 8.4]{BATYREV2} and \cite[Thm. 5.1]
{DENEF-LOESER}.)

In the present paper we give explicit formulae for the evaluation of the
function $E_{\text{str}}\left( X;u,v\right) $ for those $X$'s which are the
underlying spaces of two special series of $\mathbf{A}_{n,\ell}^{\left(
r\right) }$\textit{-singularities }(see below {(d)} for the precise
definition) by constructing an appropriate snc-resolution $\varphi:%
\widetilde{X}\longrightarrow X$, by examining the nature of the arising
exceptional prime divisors and, finally, by computing their $E$ \textit{-}%
polynomials. (In \cite{DAIS-ROCZEN} this was carried out for all
three-dimensional $\mathbf{A}$-$\mathbf{D}$-$\mathbf{E}$ singularities).

\medskip \noindent \noindent {}\textsf{(a)}\textbf{\ Log-terminal
singularities.}\thinspace ~ Let $X$ be a normal complex variety. Suppose
that $X$ is $\mathbb{Q}$-Gorenstein, i.e., that a positive integer multiple
of its canonical Weil divisor $K_{X}$ is a Cartier divisor. $X$ is said to
have at most \textit{log-terminal} (respectively, \textit{canonical}%
\thinspace /\thinspace \textit{terminal}) singularities if there exists an
\textit{snc}-desingularization $\varphi :\widetilde{X}\longrightarrow X$,
i.e., a desingularization of $X$ whose exceptional locus $\mathfrak{Ex}%
\left( \varphi \right) =\cup _{i=1}^{m}D_{i}$ consists of \textbf{s}mooth
prime divisors $D_{1},D_{2},\ldots ,D_{m}$ with only \textbf{n}ormal \textbf{%
c}rossings, such that the ``discrepancy'' w.r.t. $\varphi $ is of the form $%
K_{\widetilde{X}}-\varphi ^{\ast }\left( K_{X}\right) ={\sum\limits_{i=1}^{m}%
}\ a_{i}\ D_{i},$ with all the $a_{i}$'s $>-1$ ($\geq 0\,$/$\,>0$). These
inequalities do not depend on the particular choice of $\varphi $.

\medskip \noindent \noindent {}\textsf{(b)}\textbf{\ $E$-polynomials.}%
\thinspace ~ Deligne proved in \cite[\S 8]{DELIGNE} that the cohomology
groups $H^{i}\left( X,\mathbb{Q}\right) $ of any complex variety $X$ are
endowed with a natural \textit{mixed Hodge structure} (MHS).\textit{\ }The
same remains true if one works with cohomologies $H_{c}^{i}\left( X,\mathbb{Q%
}\right) $ with compact supports. There exist namely an increasing
weight-filtration\textit{\ }$\mathcal{W}_{\bullet }$ and a decreasing
Hodge-filtration of $H^{i}\left( X,\mathbb{Q}\right) $ (resp. $%
H_{c}^{i}\left( X,\mathbb{C}\right) )$ which induces a natural filtration $%
\mathcal{F}^{\bullet }$ on the complexification of the corresponding graded
pieces $Gr_{k}^{\mathcal{W}_{\bullet }}(H^{i}\left( X,\mathbb{Q}\right) )$
(resp. $Gr_{k}^{\mathcal{W}_{\bullet }}(H_{c}^{i}\left( X,\mathbb{Q}\right)
))$. Let
\begin{align*}
h^{p,q}(H^{i}\left( X,\mathbb{C}\right) )& :=\text{ dim}_{\mathbb{C}}Gr_{%
\mathcal{F}^{\bullet }}^{p}Gr_{p+q}^{\mathcal{W}_{\bullet }}(H^{i}\left( X,%
\mathbb{C}\right) )\text{ \smallskip\ \ } \\
\text{( resp. }h^{p,q}(H_{c}^{i}\left( X,\mathbb{C}\right) )& :=\text{ dim}_{%
\mathbb{C}}Gr_{\mathcal{F}^{\bullet }}^{p}Gr_{p+q}^{\mathcal{W}_{\bullet
}}(H_{c}^{i}\left( X,\mathbb{C}\right) )\text{ )}
\end{align*}
denote hereafter the corresponding \textit{Hodge numbers.} The so-called $E$-%
\textbf{polynomial} of $X$ is defined to be
\begin{equation*}
E\left( X;u,v\right) :=\sum_{p,q}\ e^{p,q}\left( X\right) \ u^{p}v^{q}\in
\mathbb{Z}\left[ u,v\right] ,
\end{equation*}
where $e^{p,q}\left( X\right) :=\sum_{i\geq 0}\ \left( -1\right) ^{i}\
h^{p,q}(H_{c}^{i}\left( X,\mathbb{C}\right) )$. (If we set $u=v=1,$ then $\
E\left( X;1,1\right) $ equals the usual \textit{topological} Euler
characteristic $e(X)$ of $X.$)\bigskip

\medskip \noindent \noindent {}\textsf{(c)}\textbf{\ $E_{\text{str}}$%
-functions.}\thinspace ~ To pass to string-theoretic invariants, one takes
essentially into account the ``discrepancy coefficients''.

\begin{definition}
\label{EST-DEF}Let $X$ be a normal complex variety with at most log-terminal
singularities, $\varphi :\widetilde{X}\longrightarrow X$ \ an
snc-desingularization of $X$ as in ${(\mathrm{a)}}$, $D_{1},D_{2},\ldots
,D_{m}$ the prime divisors of the exceptional locus, and $I:=\left\{
1,2,\ldots ,m\right\} $. For any subset $J\subseteq I$ define
\begin{equation*}
D_{J}:=
\begin{cases}
\widetilde{X}, & \text{if \ }J=\varnothing  \\
\, & \, \\
\bigcap_{j\in J}\ D_{j}, & \text{if \ }J\neq \varnothing
\end{cases}
\quad \text{and }\quad D_{J}^{\circ }:=D_{J}\,\mathbb{r}\bigcup_{j\in I%
\mathbb{r}J}\ D_{j}.
\end{equation*}
The algebraic function
\begin{equation}
E_{\text{str}}\left( X;u,v\right) :=\sum_{J\subseteq I}\ E\left(
D_{J}^{\circ };u,v\right) \ \prod_{j\in J}\ \frac{uv-1}{\left( uv\right)
^{a_{j}+1}-1}  \label{E-STR}
\end{equation}
(under the convention for $\prod_{j\in J}$ to be $1$, if $J=\varnothing ,$
and $E\left( \varnothing ;u,v\right) :=0$) is called the \textbf{%
string-theoretic} $E$\textbf{-function }\textrm{{(} \textit{or simply} $%
E_{str}$\textbf{-function}) \textit{of} $\ X$. }
\end{definition}

\noindent {}The major result of \cite{BATYREV1} says that:

\begin{theorem}
\label{INDEP}$E_{\text{\emph{str}}}\left( X;u,v\right) $ is independent of
the choice of the snc- desingularization $\varphi :\widetilde{X}%
\longrightarrow X$.
\end{theorem}

\begin{remark}
(i) Though the string-theoretic function $E_{\mathrm{str}}\left(
X;u,v\right) $ enjoys this particularly important invariance property, to
\textit{evaluate} it by (\ref{E-STR}) one needs not only the \textit{%
existence} of (at least one) snc-desingularization (which is guaranteed,
e.g., by Hironaka's main theorems \cite{HIRONAKA}), but also the precise
knowledge of what kind of exceptional prime divisors are available on the
corresponding smooth model, and which are their intersections. In general,
there are several ways to resolve log-terminal singularities, involving
different choices for the centers of the modifications of $X$ and,
sometimes, necessary extra normalizations, blow-ups of non-reduced
subschemes etc. For this reason, a first realistic attempt to understand the
behaviour of (\ref{E-STR}), from the computational point of view, cannot
overlook the class of \textit{absolutely isolated} singularities, i.e.,
isolated singularities resolvable by a finite sequence of (usual) blow-ups
of closed points, for which one may keep the needed details (strict
transforms after each step of the resolution procedure, snc-condition etc.)
under control.\smallskip\ \newline
{(ii)}\thinspace ~ It is also worth mentioning that the ``first summand'' in
(\ref{E-STR}), i.e., for $J=\varnothing ,$ equals
\begin{equation*}
E(\widetilde{X}\,\mathbb{r}{\textstyle\bigcup_{j=1}^{m}}\ D_{j};u,v)=E(X\,%
\mathbb{r\,}\text{Sing}(X);u,v)
\end{equation*}
(where Sing$(X)$ denotes the singular locus of $X$). This means that it can
be described exclusively by the study of topological properties of $X$
``around'' the singularities without involving any resolution data.
\end{remark}

\begin{definition}
The rational number
\begin{equation}
e_{\mathrm{str}}\left( X\right) :=\text{ }\underset{u,v\rightarrow 1}{%
\mathrm{lim}}E_{\mathrm{str}}\left( X;u,v\right) =\sum_{J\subseteq I}\
e\left( D_{J}^{\circ }\right) \ \prod_{j\in J}\ \frac{1}{a_{j}+1}\
\label{e-STR}
\end{equation}
is called the \textbf{string-theoretic Euler number} of $X$. Moreover, the
\textbf{string-theoretic index }\textrm{ind}$_{\mathrm{str}}\left( X\right) $
of $X$ is defined to be the integer
\begin{equation*}
\mathrm{ind}_{\mathrm{str}}\left( X\right) :=\mathrm{min}\left\{ \ l\in
\mathbb{Z}_{\geq 1}\mathbb{\ }\left| \ \ e_{\mathrm{str}}\left( X\right) \in
\frac{1}{l}\,\mathbb{Z\ }\right. \right\} .
\end{equation*}
\end{definition}

\begin{examples}
(i) For $\mathbb{Q}$-Gorenstein toric varieties $X$, ind$_{\text{str}}\left(
X\right) =1$, and $e_{\text{str}}\left( X\right) $ is equal to the
normalized volume of the defining fan. Moreover, for Gorenstein toric
varieties $X$, $E_{\text{str}}\left( X;u,v\right) $ is a polynomial (cf.
\cite[4.4 and 4.10]{BATYREV1}).\smallskip \newline
(ii) Normal algebraic surfaces $X$ with at most log-terminal singularities
have string-theoretic index ind$_{\text{str}}\left( X\right) =1$. There
exist, however, normal complex varieties $X$ of dimension $d\geq 3$ with at
most Gorenstein canonical singularities having ind$_{\text{str}}\left(
X\right) >1$.
\end{examples}

\noindent Batyrev's conjecture \cite[5.9]{BATYREV1}, concerning the range of
ind$_{\text{str}}\left( X\right) $, can be stated as follows:

\begin{conjecture}[Boundedness of the string-theoretic index]
\label{BAT-CONJ}Let $X$ be an $r$-dimens\-ional normal complex variety
having at most Gorenstein canonical singularities. Then ind$_{\mathrm{str}%
}\left( X\right) $ is bounded by a constant $C\left( r\right) $ depending
only on $r$.
\end{conjecture}

\noindent As it turns out (see below Remark \ref{CONJREM}), and in contrast
to initial expectations due to some classes of examples (see, e.g.,
\cite[5.1, 5.10]{BATYREV1} for the case of cones over certain smooth
projective Fano varieties), conjecture \ref{BAT-CONJ} is not true in
general. Nevertheless, the characterization of those classes of $X$'s, which
admit bounded string-theoretic indices, remains an unsolved problem.

\medskip \noindent \noindent {}\textsf{(d)}\textbf{\ The $\mathbf{A}_{n,\ell
}^{\left( r\right) }$'s.}\thinspace ~ We define the $r$-dimensional $\mathbf{%
A}_{n,\ell }^{\left( r\right) }$\textbf{-singularities} as those isolated
hypersurface singularities which have underlying spaces of the form
\begin{equation*}
X_{n,\ell }^{(r)}:=\text{Spec}\left( \mathbb{C}\left[ x_{1},\ldots ,x_{r+1}%
\right] \ /\ \text{\textbf{(}}f\text{\textbf{)}}\right) ,
\end{equation*}
where $r,n,\ell $ are integers, such that $r\geq \ell \geq 2,\ n+1\geq \ell
, $ and \fboxsep4pt
\begin{equation}
\fbox{$f\left( x_{1},\ldots ,x_{r+1}\right) :=x_{1}^{n+1}+x_{2}^{\ell
}+x_{3}^{\ell }+\cdots +x_{r+1}^{\ell }$}\,.  \label{DEFF}
\end{equation}

\noindent These are obviously singularities of \ \textit{Brieskorn-Pham type}%
. In addition, by our assumptions about $r,\ell$ and$\ n$, they are \textit{%
canonical }(see Reid \cite[Prop. 4.3, p. 297]{REID}). The notation is chosen
in this manner to remind that they include, in particular, both subseries of
usual $r$-dimensional $\mathbf{A}_{n}$-singularities ($\mathbf{A}%
_{n,2}^{\left( r\right) }$'s) and of Fermat singularities ($\mathbf{A}%
_{\ell-1,\ell}^{\left( r\right) }$'s).\bigskip

\medskip \noindent {}\noindent \textsf{(e)}\textbf{\ Some auxiliary
combinatorial functions.}\thinspace ~At first, for $p,q\in \mathbb{Z}_{\geq
0},$ let us denote Kronecker's symbol by
\begin{equation*}
\delta _{p,q}=
\begin{cases}
1, & \text{if \thinspace }p=q \\
0, & \text{if \thinspace }p\neq q.
\end{cases}
\end{equation*}
$\bullet $\thinspace ~Next, fixing $r,\ell $ and $n,$ as in \noindent
\noindent \textsf{(d)}, we set $d:=\text{lcm}\left( n+1,\ell \right) $. and
define three functions $\mathbf{a},\mathbf{b}$ and $\mathbf{c}:\mathbb{Z}%
_{\geq 0}\longrightarrow \,\mathbb{Z}_{\geq 0}$ \ by
\begin{equation}
\mathbb{Z}_{\geq 0}\ni i\longmapsto \mathbf{a}\left( i\right) ={\textstyle%
\sum\limits_{p=0}^{n-1}}\ \delta _{i,\frac{pd}{n+1}},  \label{A-FUNCTION}
\end{equation}
by the multinomial coefficients
\begin{equation}
\mathbb{Z}_{\geq 0}\ni j\longmapsto \mathbf{b}\left( j\right) =
\begin{cases}
\sum\limits_{\left( \nu _{1},\nu _{2},\ldots ,\nu _{\ell -1}\right) \in %
\mathfrak{B}_{j}}\binom{r}{\nu _{1},\nu _{2},\ldots ,\nu _{\ell -1}}, &
\text{if \thinspace }\mathfrak{B}_{j}\neq \varnothing \\
0, & \text{otherwise}
\end{cases}
\label{B-FUNCTION}
\end{equation}
(with $\binom{r}{\nu _{1},\nu _{2},\ldots ,\nu _{\ell -1}}:=\frac{r!}{\nu
_{1}!\,\nu _{2}!\,\ldots \,\nu _{\ell -1}!}$), and by the convolutional
formula
\begin{equation}
\mathbb{Z}_{\geq 0}\ni k\longmapsto \mathbf{c}\left( k\right) =\sum_{\left(
i,j\right) \in \mathfrak{C}_{k}}\,\mathbf{a}\left( i\right) \,\mathbf{b}%
\left( j\right) ,  \label{CONVOLUTION}
\end{equation}
where for each $j\in \mathbb{Z}_{\geq 0},$

{\small
\begin{equation*}
\mathfrak{B}_{j}:=\left\{ \!\left(
\nu_{1},\nu_{2},\ldots,\nu_{\ell-1}\right) \in\left( \mathbb{Z}%
_{\geq0}\right) ^{\ell-1}\ \!\left|
\begin{array}{l}
\nu_{1}+\nu_{2}+\cdots+\nu_{\ell-1}=r\text{\ \ } \\[3pt]
\text{and (whenever }\ell\geq3) \\[3pt]
d\left( \nu_{2}+2\nu_{3}+\cdots+\left( \ell-2\right) \nu_{\ell-1}\right)
=j\,\ell
\end{array}
\right. \!\right\} ,
\end{equation*}
} \smallskip \noindent and for each $k\in\mathbb{Z}_{\geq0},$ $\mathfrak{C}%
_{k}:=\left\{ \left( i,j\right) \in\left( \mathbb{Z}_{\geq0}\right) ^{2}\
\left| \ i+j=k\right. \right\}$.

\noindent{}$\bullet$\,~Finally, for any four-tuple $\left( \kappa,\lambda
,\nu,\xi\right) \in\left( \mathbb{Z}_{\geq0}\right) ^{4}$ with $\kappa
\geq\lambda,$ let us recall the definition of the \textit{non-central
Eulerian numbers }$\mathfrak{S}\left( \kappa,\lambda\ \left| \
\nu,\xi\right. \right) $ \textit{of generalized factorials }(with
translation summand $\xi$). These are the coefficients which occur in the
development of
\begin{equation*}
\textstyle{\binom{\nu\cdot t+\xi}{\kappa}=\sum\limits_{\lambda=0}^{\kappa}\ %
\mathfrak{S}\left( \kappa,\lambda\ \left| \ \nu,\xi\right. \right) \ \binom{%
t+\kappa-\lambda }{\kappa}}
\end{equation*}
and satisfy the recurrence relation
\begin{align*}
\left( \kappa+1\right) \,\mathfrak{S}\left( \kappa+1,\lambda\ \left| \ \nu
,\xi\right. \right) ={}&\left( \nu\lambda-\kappa+\xi\right) \,\mathfrak{S}%
\left(\kappa,\lambda\ \left| \ \nu,\xi\right. \right) \\
&{}+\left( \nu\left( \kappa-\lambda+1\right) +\kappa-\xi\right) \,%
\mathfrak{S}\left( \kappa,\lambda-1\ \left| \ \nu,\xi\right. \right)
\end{align*}
with initial conditions $\mathfrak{S}\left( 0,0\ \left| \ \nu,\xi\right.
\right) =1$ and $\mathfrak{S}\left( \kappa,0\ \left| \ \nu,\xi\right.
\right) =\binom{\xi}{\kappa}$. In fact, it can be shown that
\begin{equation*}
\textstyle{\mathfrak{S}\left( \kappa,\lambda\ \left| \ \nu,\xi\right.
\right) =\sum\limits _{j=0}^{\lambda}\,\left( -1\right) ^{j}}\,\textstyle{%
\binom{\kappa+1}{j}}\,\textstyle{\binom {\nu\left( \lambda-j\right)
+\xi}{\kappa}}.
\end{equation*}

\medskip \noindent \noindent {}\textsf{(f)}\textbf{\ Main results.}%
\thinspace ~ We can now state the main results.

\begin{proposition}
\label{PROP-FS}Let $X=X_{n,\ell }^{(r)}$ be the underlying space of the $%
\mathbf{A}_{n,\ell }^{\left( r\right) }$- singularity. Then the $E$
-polynomial $E\left( X\mathbb{r}\{\mathbf{0}\};u,v\right) $ equals
{\footnotesize
\begin{equation}
\left( uv-1\right) \left[ 1+\left( uv\right) ^{r-1}+{\textstyle%
\sum\limits_{p=1}^{r-2}}\big(\left( uv\right) ^{p}+\left( -1\right) ^{r}%
\mathbf{c}(d(p+\tfrac{n}{n+1}-\tfrac{r}{\ell }))u^{p}\,v^{r-p-1}\big)\right]
\label{FIRST SUMMAND}
\end{equation}
} $($with the $\mathbf{c}$-function as defined in \emph{(\ref{CONVOLUTION})).%
}
\end{proposition}

\noindent {}Formula (\ref{FIRST SUMMAND}) provides the ``first summand'' of
the $E_{\text{str}}$-function of $X$. On the other hand, if $\ell $ divides
either $n$ or $n+1,$ $\mathbf{A}_{n,\ell }^{\left( r\right) }$'s are
absolutely isolated (see below Proposition \ref{DESINGULARIZATION}), and the
$E_{\text{str}}$-function of $X_{n,\ell }^{(r)}$'s is computed as follows:

%==============MAIN THEOREM======================

\begin{theorem}
\label{MAIN} If the integer $\ell $ divides either $n$ or $n+1,$ then the $%
E_{\text{\emph{str}}}$-function of $X=X_{n,\ell }^{(r)}$ is given by the
formula{\footnotesize {\scriptsize \thinspace
\begin{equation*}
\begin{tabular}{|c|c|}
\hline
{\small {$\mathbf{Cases}$}} & $
\begin{array}{c}
\, \\
\text{{\small {$E_{\emph{str}}\left( X;u,v\right) $}}}\allowbreak  \\
\,
\end{array}
$ \\ \hline\hline
$\ell \left| n+1\right. $ & $\!\!\!\!
\begin{array}{l}
\  \\
E\left( X\mathbb{r}\{\mathbf{0}\};u,v\right)  \\
\! \\
+\left( uv-1\right) \,\,\biggl(\tfrac{uv}{\left( uv\right) ^{r-\ell +1}-1}+{%
\textstyle}\sum\limits_{i=2}^{m-1}\tfrac{uv-1}{\left( uv\right) ^{i(r-\ell
)+1}-1}-\tfrac{uv-1}{\left( uv\right) ^{m(r-\ell )+1}-1} \\
\  \\
{\textstyle}\sum\limits_{i=1}^{m-1}\tfrac{uv-1\ }{\left( \left( uv\right)
^{i(r-\ell )+1}-1\right) \left( \left( uv\right) ^{(i+1)(r-\ell
)+1}-1\right) }\biggl) \\
\  \\
\times \left[ {\textstyle\sum\limits_{p=0}^{r-2}}\,u^{p}\,(v^{p}+\left(
-1\right) ^{r-2}\mathfrak{S}\left( r-1,p+1\ \left| \ \ell -1,p\right.
\right) \,v^{r-2-p})\right]  \\
\, \\
+\tfrac{\left( uv-1\right) }{\left( uv\right) ^{m(r-\ell )+1}-1}\left[ {%
\textstyle\sum\limits_{p=0}^{r-1}}\,u^{p}\,(v^{p}+\left( -1\right) ^{r-1}%
\mathfrak{S}\left( r,p+1\ \left| \ \ell -1,p\right. \right) \,v^{r-1-p})%
\right]  \\
\  \\
\text{\emph{\lbrack In the special case in which }}\ell =n+1\text{\emph{,
one has to} delete\emph{\ the 2nd summand}.\emph{]}} \\
\
\end{array}
$ \\ \hline
$\ell \left| n\right. $ & $\!\!\!\!
\begin{array}{c}
\begin{array}{l}
\, \\
E\left( X\mathbb{r}\{\mathbf{0}\};u,v\right) +\tfrac{uv)^{r}-1}{\left(
uv\right) ^{\left( m-1\right) \,\ell \,\left( r-\ell \right) +\,r}-1} \\
\! \\
\begin{array}{l}
+\left( uv-1\right) \,\,\biggl(\tfrac{uv}{\left( uv\right) ^{r-\ell +1}-1}+{%
\textstyle}\sum\limits_{i=2}^{m-1}\tfrac{uv-1}{\left( uv\right) ^{i(r-\ell
)+1}-1} \\
\, \\
-\tfrac{uv-1}{\left( uv\right) ^{\left( m-1\right) \,\ell \,\left( r-\ell
\right) +\,r}-1}+{\textstyle}\sum\limits_{i=1}^{m-2}\tfrac{uv-1}{\left(
\left( uv\right) ^{i(r-\ell )+1}-1\right) \left( \left( uv\right)
^{(i+1)(r-\ell )+1}-1\right) } \\
\, \\
+\tfrac{uv-1}{\left( \left( uv\right) ^{(m-1)(r-\ell )+1}-1\right) \left(
\left( uv\right) ^{\left( m-1\right) \,\ell \,\left( r-\ell \right)
+\,r}-1\right) }\biggl) \\
\  \\
\times \left[ {\textstyle\sum\limits_{p=0}^{r-2}}\,u^{p}\,(v^{p}+\left(
-1\right) ^{r-2}\mathfrak{S}\left( r-1,p+1\ \left| \ \ell -1,p\right.
\right) \,v^{r-2-p})\right]
\end{array}
\
\end{array}
\  \\
\,
\end{array}
$ \\ \hline
\end{tabular}
\medskip
\end{equation*}
}\smallskip } \noindent In particular, for the string-theoretic Euler number
we obtain\emph{:}\smallskip {\scriptsize {\small
\begin{equation*}
\begin{tabular}{|c|c|}
\hline
$\mathbf{Cases}$ & $e_{\text{\emph{str}}}\left( X\right) $ \\ \hline\hline
$\ell \left| n+1\right. $ & $\!
\begin{array}{c}
\  \\
\begin{array}{l}
\frac{m-1}{m\left( r-\ell \right) +1}\allowbreak \left[ \tfrac{1}{\ell }%
\left( \left( 1-\ell \right) ^{r}-1\right) +r\right]  \\
\, \\
+\tfrac{1}{m\left( r-\ell \right) +1}\left[ \tfrac{1}{\ell }(\left( 1-\ell
\right) ^{r+1}-1)+r+1\right]
\end{array}
\\
\
\end{array}
\!$ \\ \hline
$\ell \left| n\right. $ & $
\begin{array}{c}
\!
\begin{array}{l}
\  \\
\frac{r}{\left( m-1\right) \,\ell \,\left( r-\ell \right) +\,r}+\allowbreak
\frac{\left( m-1\right) \ell }{\left( r-\ell \right) \left( m-1\right) \ell
+r}\left[ \tfrac{1}{\ell }\left( \left( 1-\ell \right) ^{r}-1\right) +r%
\right]  \\
\,
\end{array}
\!
\end{array}
$ \\ \hline
\end{tabular}
\smallskip
\end{equation*}
} }The above number $m$ is defined to be
\begin{equation}
m:=\left\{
\begin{array}{ll}
\frac{n+1}{\ell }, & \text{\emph{if} \thinspace }n+1\equiv 0\left( \mathrm{{%
mod}\text{ }\ell }\right)  \\
\  &  \\
\frac{n}{\ell }+1, & \text{\emph{if} \thinspace }n\equiv 0\left( \mathrm{{mod%
}\text{ }\ell }\right)
\end{array}
\right.   \label{DEF-M}
\end{equation}
\bigskip
\end{theorem}

\begin{remark}
\label{CONJREM} Counterexamples to conjecture \ref{BAT-CONJ} occur already
for $\ell =2,$ as we have:
\begin{equation*}
e_{\text{str}}\left( X_{n,2}^{\left( r\right) }\right) =\left\{
\begin{array}{ll}
\frac{m(r-1)+2}{m(r-2)+1}=\frac{n\left( r-1\right) +r+3}{n(r-2)+r}, & \text{%
if both }n\text{ and }r\text{ are odd } \\
\  & \  \\
\frac{m\,r}{m\left( r-2\right) +1}=\frac{r(n+1)}{(r-2)(n+1)+2}, & \text{if }n%
\text{ is odd and }r\text{ even} \\
\  & \  \\
\frac{2(m-1)(r-1)+r}{2(m-1)(r-2)+r}=\frac{(r-1)n+r}{(r-2)n+r}, & \text{if }n%
\text{ is even and }r\text{ odd} \\
\  & \  \\
\frac{(2m-1)r}{2(m-1)(r-2)+r}=\frac{r(n+1)}{(r-2)n+r}, & \text{if both }n%
\text{ and }r\text{ are even}
\end{array}
\right.
\end{equation*}
For instance, in dimension $r=3$, we obtain:
\begin{equation*}
\underset{n\rightarrow \infty ,\text{ }n\text{ even}}{\text{lim}}\text{ ind}%
_{\text{str}}\left( X_{n,2}^{\left( 3\right) }\right) =\infty .
\end{equation*}
On the other hand, for \textit{all} odd $n$'s, $e_{\text{str}%
}(X_{n,2}^{\left( 3\right) })=2$ and ind$_{\text{str}}(X_{n,2}^{\left(
3\right) })=1.$
\end{remark}

\section{On the MHS of the cohomology groups \newline
of links}

\noindent{}At first, we shall exploit the fact that $\mathbf{A}_{n,\ell
}^{\left( r\right) }$'s are quasihomogeneous singularities, and show that
Proposition \ref{PROP-FS} is a byproduct of a more general result concerning
isolated singularities of this sort (see \ref{BIG Prop}).\medskip

\noindent{}\noindent\textsf{(a) }\textbf{Links and Milnor fibers. }Let $%
\left( W,\mathbf{0}\right) \subseteq\left( \mathbb{C}^{N},\mathbf{0}\right) $
be the germ of a complex analytic set $W$ having pure dimension $r+1$ and
the origin as isolated singularity. Assume that $f:W\rightarrow \mathbb{C}$
is a holomorphic function, such that $f\left| _{W\mathbb{r}\{\mathbf{0}%
\}}\right. $ is non-singular. Obviously, $X:=f^{-1}\left( \mathbf{0}\right) $
is a complex analytic subset of $\mathbb{C}^{N}$ of pure dimension $r$ with
the origin as isolated singularity. Let $L:=L(X,\mathbf{0}):=\mathbb{S}%
_{\varepsilon}\cap X$ denote its \textit{link}, where $\mathbb{S}%
_{\varepsilon}:=\{\mathbf{z}\in\mathbb{C}^{N}\ \left| \ \left\| \mathbf{z}%
\right\| =\varepsilon\right. \}$, $0<\varepsilon\ll1.$ $L$ is a
differentiable, compact, oriented manifold of dimension $2r-1$, and there
are isomorphisms:
\begin{equation}
H^{i+1}\left( X,X\mathbb{r}\left\{ \mathbf{0}\right\} ,\mathbb{Q}\right)
\cong H^{i}\left( X\mathbb{r}\left\{ \mathbf{0}\right\} ,\mathbb{Q}\right)
\cong H^{i}\left( L,\mathbb{Q}\right) \ .  \label{ISOL-X}
\end{equation}
If $\mathbb{B}_{\varepsilon^{\prime}}$ is the open ball with $\mathbf{0}$ as
its center and $\varepsilon^{\prime}$ as its radius, where $\varepsilon
<\varepsilon^{\prime}\ll1$, it is known that the map
\begin{equation*}
f\left| _{\mathbb{B}_{\varepsilon^{\prime}}\cap f^{-1}(\mathbb{D}_{\alpha
}^{\ast})}\right. :\mathbb{B}_{\varepsilon^{\prime}}\cap f^{-1}(\mathbb{D}%
_{\alpha}^{\ast})\longrightarrow\mathbb{D}_{\alpha}^{\ast}
\end{equation*}
determines a differentiable fibre bundle, where $\mathbb{D}_{\alpha}^{\ast
}:=\{t\in\mathbb{C}\ \left| \ 0<\left| t\right| <\alpha\right. \}$ is a
small punctured disc in $\mathbb{C}$ with $0<\alpha<\varepsilon.$ Let $%
F=F_{t}$ be the corresponding fiber, the so-called (\textit{open}) \textit{%
Milnor fiber}. The study of the relation between the MH-structures of the
cohomology groups of $L$ and $F$ relies on certain corollaries of a theorem
of Steenbrink \cite[(2.3)]{STEENBRINK2} and Hamm \cite[Thm. 1.6.1]{HAMM3}.
(The coefficients of the cohomology groups are always taken from $\mathbb{C}%
. $)

\begin{theorem}[Steenbrink-Hamm]
\label{SH-THM}For all $i,$ there exists an exact MHS-sequence\emph{:}
\begin{equation*}
\cdots \longrightarrow H^{i-1}\left( L\right) \longrightarrow
H_{c}^{i}\left( F\right) \longrightarrow H^{i}\left( F\right)
\longrightarrow H^{i}\left( L\right) \longrightarrow \cdots
\end{equation*}
\end{theorem}

\begin{corollary}
\label{COROL1}We have the following exact sequence and isomorphisms of MHS%
\emph{:}\medskip \newline
{\small $
\begin{array}{ll}
\emph{(i)} & 0\!\rightarrow \!H^{r-1}\left( F\right) \!\rightarrow
\!H^{r-1}\left( L\right) \!\rightarrow \!H_{c}^{r}\left( F\right)
\!\rightarrow \!H^{r}\left( F\right) \!\rightarrow \!H^{r}(L)\!\rightarrow
\!H_{c}^{r+1}\left( F\right) \!\rightarrow 0\ \smallskip  \\
\emph{(ii)} & \smallskip H^{i}\left( L\right) \cong H^{i}\left( F\right) ,%
\text{ \ for all }i<r-1, \\
\emph{(iii)} & \smallskip H^{i-1}\left( L\right) \cong H_{c}^{i}\left(
F\right) ,\text{ \ for \ all \ }i>r-1.
\end{array}
$\emph{\ } }
\end{corollary}

\noindent{}\textsc{Proof}\textit{. }Since $\mathbb{B}_{\varepsilon^{\prime}}$
is a complex Stein manifold, $F$ is a complex Stein manifold too. Hence, $F$
has the homotopy type of a CW-complex of real dimension $r$ (see \cite{HAMM2}%
), which means that $H^{i}(F)\cong H_{c}^{2r-i}(F)=0$ for all $i\geq r+1.$
The exactness in (i) and the existence of MHS-isomorphisms (ii) and (iii)
follow from the long exact sequence of Theorem \ref{SH-THM}, combined with
the vanishing of these cohomology groups.\hfill$\square$

\begin{corollary}
\label{COROL2}For all $p,q,$ the Hodge numbers of the two ``middle''
cohomology groups of $F$ satisfy the equalities\textbf{:}
\begin{align*}
h^{p,q}(H^{r}(F))&
=h^{p,q}(H^{r-1}(F))+h^{p,q}(H_{c}^{r}(F))-h^{p,q}(H_{c}^{r+1}(F))+ \\
& +h^{p,q}(H^{r}(L))-h^{p,q}(H^{r-1}(L))\smallskip  \\
& =h^{p,q}(H^{r-1}(F))+h^{r-p,r-q}(H^{r}(F))-h^{r-p,r-q}(H^{r-1}(F))+ \\
& +h^{p,q}(H^{r}(L))-h^{p,q}(H^{r-1}(L))
\end{align*}
\end{corollary}

\noindent{}{}\textsc{Proof}\textit{. }The first equality is obvious by \ref
{COROL1} (i), and the second one follows from Poincar\'{e} duality.\hfill $%
\square$

\begin{proposition}
\label{COROL3}If $N=r+1,W=\mathbb{C}^{r+1}$ and $(X,\mathbf{0})$ is a purely
$r$-dimensional isolated hypersurface singularity, with $r\geq 2,$ then the
only ``non-trivial'' Hodge numbers of the cohomology groups of its link $%
L=L(X,\mathbf{0})$ are
\begin{equation*}
h^{p,q}(H^{r-1}(L))=h^{r-p,r-q}(H^{r}(L)),\text{ with }p+q\leq r-1,
\end{equation*}
as we have\emph{:}\medskip \newline
$
\begin{array}{ll}
\emph{(i)} & h^{p,q}(H^{i}(L))=0,\text{ for all }p,q\text{ whenever }i\notin
\left\{ 0,r-1,r,2r-1\right\} .\smallskip  \\
\emph{(ii)} & h^{p,q}(H^{0}(L))=1,\text{ for }p\smallskip =q=0,\text{ and}=0,%
\text{ otherwise.\newline
\ \smallskip } \\
\emph{(iii)} & h^{p,q}(H^{2r-1}(L))=1,\text{ for }p\smallskip =q=r,\text{
and }=0,\text{ otherwise.\thinspace \smallskip \smallskip } \\
\emph{(iv)} &
\begin{array}[t]{l}
h^{p,q}(H^{r-1}(L))=h^{r-p,r-q}(H^{r}(L)),\text{ \ for all }p,q,\smallskip
\\
\text{and equals }0\text{ whenever }p+q>r-1.
\end{array}
\end{array}
$ \newline
\end{proposition}

\noindent{}\textsc{Proof}\textit{. }$L$ is $\left( r-2\right) $-connected
(cf. \cite[Thm. 5.2]{MILNOR}), and the local Lefschetz Theorem gives (i),
(ii) and (iii) because $H^{i}(L)=0$ \ for all indices $i\notin\left\{
0,r-1,r,2r-1\right\} $ and $H^{0}(L)\cong H^{2r-1}(L)\cong\mathbb{C}.$ For
(iv) use Poincar\'{e} duality and the fact, that the natural MHS on $%
H^{i}\left( L\right) $ has weights $Gr_{j}^{\mathcal{W}_{\bullet}}(H^{i}%
\left( L\right) )=0$ for $j>i$ (by the Semipurity Theorem, cf. \cite[Cor.
(1.12), p. 518]{STEENBRINK2}).\hfill$\square\bigskip$

\noindent\textsf{(b) }\textbf{Quasihomogeneous isolated singularities. }A
polynomial
\begin{equation*}
f\in\mathbb{C}\left[ x_{1},x_{2},\ldots,x_{r+1}\right]
\end{equation*}
is \textit{quasihomogeneous} of degree $d$ with respect to the \textit{%
weights}\
\begin{equation*}
\mathbf{w}=\left( w_{1},\ldots,w_{r+1}\right) \in\left( \mathbb{Z}_{\geq
1}\right) ^{r+1}
\end{equation*}
if
\begin{equation*}
f\left( \lambda^{w_{1}}x_{1},\ldots,\lambda^{w_{r+1}}x_{r+1}\right)
=\lambda^{d}\ f\left( x_{1},\ldots,x_{r+1}\right) ,\ \ \ \forall \lambda,\ \
\ \lambda\in\mathbb{C}^{\ast}.
\end{equation*}
Hereafter we consider such an $f,$ assume that $r\geq2$ and that
\begin{equation*}
X_{f}:=\left\{ \left( x_{1},\ldots,x_{r+1}\right) \in\mathbb{C}^{r+1}\
\left| \ \right. f\left( x_{1},\ldots,x_{r+1}\right) =0\right\}
\end{equation*}
has no other singularities than\textbf{\ }$\mathbf{0}\in\mathbb{C}^{r+1}.$
Note that the Milnor algebra
\begin{equation*}
M\left( f\right) :=\mathbb{C}\left[ x_{1},x_{2},\ldots,x_{r+1}\right] \ /\
\left( \frac{\partial f}{\partial x_{1}},\ldots,\frac{\partial f}{\partial
x_{r+1}}\right)
\end{equation*}
associated to $f$ is a graded $\mathbb{C}$-algebra of finite type \ (with deg%
$\left( x_{i}\right) =w_{i}$, $\ i=1,...,r+1$) whose Poincar\'{e} series $%
P_{M\left( f\right) }\left( t\right) $ equals
\begin{equation}
P_{M\left( f\right) }\left( t\right) ={\textstyle\sum_{k\geq0}} \,\dim_{%
\mathbb{C}}\left( M\left( f\right) _{k}\right) \,t^{k}=\tfrac{\left(
1-t^{d-w_{1}}\right) \left( 1-t^{d-w_{2}}\right) \cdots\left(
1-t^{d-w_{r+1}}\right) }{\left( 1-t^{w_{1}}\right) \left( 1-t^{w_{2}}\right)
\cdots\left( 1-t^{w_{r+1}}\right) }  \label{POIN-POL}
\end{equation}
(cf. \cite[(7.27), p. 112 ]{DIMCA1}). Next, we define the \textit{quasismooth%
} weighted projective hypersurfaces
\begin{equation*}
Z=\left\{ \left[ x_{0}:x_{1}:\ldots:x_{r+1}\right] \in\mathbb{P}_{\mathbb{C}%
}^{r+1}\left( 1,\mathbf{w}\right) \ \left| \ \right. \overline{f}\left(
x_{0},\ldots,x_{r+1}\right) =0\right\}
\end{equation*}
where $\overline{f}\left( x_{0},\ldots,x_{r+1}\right) :=x_{0}^{d}-f\left(
x_{1},\ldots,x_{r+1}\right) $, and
\begin{align*}
Z_{\infty} & =\left\{ \left[ x_{0}:x_{1}:\ldots:x_{r+1}\right] \in Z\ \left|
\ \right. x_{0}=0\right\} \\
& \cong\left\{ \left[ x_{1}:\ldots:x_{r+1}\right] \in\mathbb{P}_{\mathbb{C}%
}^{r}\left( \mathbf{w}\right) \ \left| \ \right. f\left(
x_{1},\ldots,x_{r+1}\right) =0\right\} .
\end{align*}
We have
\begin{equation}
M\left( \overline{f}\right) =M\left( f\right) \otimes\mathbb{C}\left[ x_{0}%
\right] \,/\,\left( x_{0}^{d-1}\right)  \label{MIL-F-FDACH}
\end{equation}
and the map $\left( x_{1},\ldots,x_{r+1}\right) \longmapsto$ $\left[
1:x_{1}:\ldots:x_{r+1}\right] $ induces a diffeomorphism between
\begin{equation*}
F=\left\{ \left( x_{1},\ldots,x_{r+1}\right) \in\mathbb{C}^{r+1}\ \left| \
\right. f\left( x_{1},\ldots,x_{r+1}\right) =1\right\}
\end{equation*}
and the complement $Z\mathbb{r}Z_{\infty},$ where this $F$ is diffeomorphic
to the (usual) Milnor fiber of the singularity $(X_{f},\mathbf{0})$ (see
\cite[(1.13), p. 72]{DIMCA2}). Moreover, $F$ has the homotopy type of a
bouquet of $\mu\left( f\right) $ $r$-spheres, with
\begin{equation}
\mu\left( f\right) =\text{ }\underset{t\rightarrow1}{\lim}P_{M\left(
f\right) }\left( t\right) ={\textstyle\prod\limits_{i=1}^{r+1}} \left(
\tfrac{d}{w_{i}}-1\right)  \label{MILNOR-NUMBER}
\end{equation}
denoting the corresponding \textit{Milnor number}. The \textit{primitive}
cohomology groups of $Z_{\infty}$ are defined by the exact sequence
\begin{equation*}
0\longrightarrow H^{r-1}(\mathbb{P}_{\mathbb{C}}^{r}\left( \mathbf{w}\right)
,\mathbb{C})\longrightarrow H^{r-1}(Z_{\infty},\mathbb{C})\longrightarrow H_{%
\text{prim}}^{r-1}\left( Z_{\infty},\mathbb{C}\right) \longrightarrow 0\ .
\end{equation*}
Since both $\mathbb{P}_{\mathbb{C}}^{r}\left( \mathbf{w}\right) $ and $%
Z_{\infty}$ are orbifolds, they are equipped with pure Hodge structure, and
therefore both $H^{r-1}(Z_{\infty},\mathbb{C})$ $\ $and $H_{\text{prim}%
}^{r-1}\left( Z_{\infty},\mathbb{C}\right) $ decompose, say as
\begin{equation*}
H^{r-1}\left( Z_{\infty},\mathbb{C}\right)
=\!\bigoplus_{p+q=r-1}H^{p,q}\left( Z_{\infty}\right) ,H_{\text{prim}%
}^{r-1}\left( Z_{\infty },\mathbb{C}\right) =\!\bigoplus_{p+q=r-1}H_{\text{%
prim}}^{p,q}\left( Z_{\infty}\right) ,
\end{equation*}
(The same is also valid for $H_{[\text{prim]}}^{r}\left( Z,\mathbb{C}\right)
$).

\begin{lemma}
\label{MIL-LEMMA}For the Milnor fiber $F$ of $(X_{f},\mathbf{0})$ we
have\medskip\ \newline
$
\begin{array}{ll}
\text{\emph{(i)}} & h^{p,q}(H^{0}(F,\mathbb{C}))=1,\text{ for }p\smallskip
=q=0,\text{ and}=0,\text{ otherwise}.\smallskip  \\
\text{\emph{(ii)}} & H^{i}\left( F,\mathbb{C}\right) =0,\text{ for all }%
i\notin \left\{ 0,r\right\} \text{.}\smallskip  \\
\text{\emph{(iii)}} & h^{p,q}\left( H^{r}\left( F,\mathbb{C}\right) \right)
=0,\,\text{\ for }p+q\notin \left\{ r,r+1\right\} .\ \smallskip  \\
\text{\emph{(iv)}} & h^{p,r-p}\left( H^{r}\left( F,\mathbb{C}\right) \right)
=\smallskip h_{\text{\emph{prim}}}^{p,r-p}\left( Z\right) =h^{p,r-p}\left(
Z\right) -\delta _{p,r-p},\,\text{\ for }0\leq p\leq r. \\
\text{\emph{(v)}} &
\begin{array}[t]{l}
h^{p,r+1-p}\left( H^{r}\left( F,\mathbb{C}\right) \right) =h_{\text{\emph{%
prim}}}^{p-1,r-p}\left( Z_{\infty }\right) = \\
=h^{p-1,r-p}\left( Z_{\infty }\right) -\delta _{p-1,r-p},\,\text{\ for }%
1\leq p\leq r.
\end{array}
\end{array}
$ \newline
\end{lemma}

\noindent\textsc{Proof}. (i) This follows from \ref{COROL1} (ii) and \ref
{COROL3} (ii).\smallskip

\noindent{}(ii)-(v). At first note that $H^{p,q}(\mathbb{P}_{\mathbb{C}%
}^{r}\left( \mathbf{w}\right) )$ (resp., $H^{p,q}(\mathbb{P}_{\mathbb{C}%
}^{r+1}\left( 1,\mathbf{w}\right) )$) is $\cong\mathbb{C},$ whenever $p=q,$
and $=0,$ otherwise. As Steenbrink points out in \cite[p. 216]{STEENBRINK1},
there is an exact MHS-sequence of Gysin-type:{\small
\begin{equation*}
\cdots\!\rightarrow H^{i}\left( Z,\mathbb{C}\right) \rightarrow H^{i}\left( Z%
\mathbb{r}Z_{\infty},\mathbb{C}\right) \rightarrow H^{i-1}\left( Z_{\infty },%
\mathbb{C}\right) (-1)\overset{\theta}{\longrightarrow}H^{i+1}\left( Z,%
\mathbb{C}\right) \rightarrow\!\cdots
\end{equation*}
} By the Weak Lefschetz Theorem \cite[4.2.2]{DOLGACHEV}, the homomorphism $\
$%
\begin{equation*}
H^{p,q}\left( Z_{\infty}\right) \!\longrightarrow H^{p,q+1}\left( \mathbb{P}%
_{\mathbb{C}}^{r}\left( \mathbf{w}\right) \right) \text{ [resp., }%
H^{p,q}\left( Z\right) \longrightarrow H^{p,q+1}\left( \mathbb{P}_{\mathbb{C}%
}^{r+1}\left( 1,\mathbf{w}\right) \right) \text{]}
\end{equation*}
is an isomorphism for $p+q>r-1$ (resp., $p+q>r$) and an epimorphism for $%
p+q=r-1$ (resp., $p+q=r$).$\ $Thus, $\theta$ is an isomorphism for all $%
i\notin\left\{ 0,r\right\} ,$ proving (ii). Moreover, since
\begin{equation*}
\mathcal{W}_{j}\left( H^{r}\left( F,\mathbb{C}\right) \right) =\left\{
\begin{array}{ll}
0, & \text{if }j<r \\
H^{r}\left( F,\mathbb{C}\right) , & \text{if }j>r
\end{array}
\right.
\end{equation*}
i.e., $Gr_{j}^{\mathcal{W}_{\bullet}}\left( H^{r}\left( F,\mathbb{C}\right)
\right) =0,$ for $j\notin\{r,r+1\},$ (cf. \cite[\S8.2]{DELIGNE}), (iii) is
obvious, and the above exact MHS-sequence gives the isomorphisms
\begin{align*}
\mathcal{W}_{r}\left( H^{r}\left( F,\mathbb{C}\right) \right) & =\text{Im}%
\left( H^{r}\left( Z,\mathbb{C}\right) \longrightarrow H^{r}\left( Z\mathbb{r%
}Z_{\infty},\mathbb{C}\right) \right) \\
& \cong\text{CoKer}\left( H^{r-2}\left( Z_{\infty},\mathbb{C}\right)
(-1)\longrightarrow H^{r}\left( Z,\mathbb{C}\right) \right) \\
& \cong\text{CoKer}\left( H^{r}(\mathbb{P}_{\mathbb{C}}^{r+1}\left( 1,%
\mathbf{w}\right) )\longrightarrow H^{r}\left( Z,\mathbb{C}\right) \right)
=H_{\text{prim}}^{r}\left( Z,\mathbb{C}\right)
\end{align*}
and
\begin{equation*}
\begin{array}{l}
Gr_{r+1}^{\mathcal{W}_{\bullet}}\left( H^{r}\left( F,\mathbb{C}\right)
\right) \\
\  \\
=H^{r}\left( Z\mathbb{r}Z_{\infty},\mathbb{C}\right) \ /\ \text{Ker}\left(
H^{r}\left( Z\mathbb{r}Z_{\infty},\mathbb{C}\right) \longrightarrow
H^{r-1}\left( Z_{\infty},\mathbb{C}\right) (-1)\right) \\
\  \\
\cong\text{Ker}\left( H^{r-1}\left( Z_{\infty},\mathbb{C}\right)
(-1)\longrightarrow H^{r+1}\left( Z,\mathbb{C}\right) \right) \\
\  \\
\cong\text{CoKer}\left( H^{r-1}(\mathbb{P}_{\mathbb{C}}^{r}\left( \mathbf{w}%
\right) ,\mathbb{C})(-1)\longrightarrow H^{r-1}\left( Z_{\infty },\mathbb{C}%
\right) (-1)\right) \\
\  \\
=H_{\text{prim}}^{r-1}\left( Z_{\infty},\mathbb{C}\right) \left( -1\right) ,
\end{array}
\end{equation*}
respectively, proving (iv) and (v).\hfill$\square$

\begin{theorem}[Griffiths-Steenbrink]
If $(X_{f},\mathbf{0})$ is an $r$-dimensional isolated quasihomogeneous
hypersurface singularity of degree $d$ w.r.t. the weights $w_{1},\ldots
,w_{r+1},$ then
\begin{equation*}
H_{\emph{prim}}^{p-1,r-p}\left( Z_{\infty }\right) \cong M\left( f\right)
_{pd-\left( w_{1}+\ldots +w_{r+1}\right) }\ .
\end{equation*}
Hence,
\begin{equation}
\left\{
\begin{array}{l}
h^{p,r-p}\left( H^{r}\left( F,\mathbb{C}\right) \right)
=\sum\limits_{i=1}^{d-1}\dim _{\mathbb{C}}\left( M\left( f\right)
_{pd-\left( w_{1}+\ldots +w_{r+1}\right) +i}\right)  \\
\  \\
h^{p+1,r-p}\left( H^{r}\left( F,\mathbb{C}\right) \right) =\dim _{\mathbb{C}%
}\left( M\left( f\right) _{(p+1)d-\left( w_{1}+\ldots +w_{r+1}\right)
}\right)
\end{array}
\right.   \label{HPQ-F}
\end{equation}
\end{theorem}

\noindent\textsc{Proof}. Extending Griffiths' results \cite{GRIFFITHS} to
the case of weighted homogeneous hypersurfaces, the global sections of the
sheaves
\begin{equation*}
\Omega_{\mathbb{P}_{\mathbb{C}}^{r}\left( \mathbf{w}\right) }^{p}(Z_{\infty
})=\Omega_{\mathbb{P}_{\mathbb{C}}^{r}\left( \mathbf{w}\right) }^{p}\otimes%
\mathcal{O}_{\mathbb{P}_{\mathbb{C}}^{r}\left( \mathbf{w}\right)
}(Z_{\infty}),
\end{equation*}
as well as the graded pieces of middle cohomology of $\mathbb{P}_{\mathbb{C}%
}^{r}\left( \mathbf{w}\right) \mathbb{r}Z_{\infty},$ are described by means
of special auxiliary differential forms with poles along $Z_{\infty}.$ In
particular,
\begin{equation*}
\begin{array}{l}
H^{0}(\mathbb{P}_{\mathbb{C}}^{r}\left( \mathbf{w}\right) ,\Omega _{\mathbb{P%
}_{\mathbb{C}}^{r}\left( \mathbf{w}\right) }^{r}(Z_{\infty }))=\left\{
\left. \dfrac{g\cdot\Omega_{0}}{f}\right| g\in\mathbb{C}%
[x_{1},..,x_{r+1}]_{d-\left( w_{1}+\ldots+w_{r+1}\right) }\right\} ,
\end{array}
\end{equation*}
where
\begin{equation*}
\begin{array}{c}
\Omega_{0}:=\sum\limits_{i=1}^{r+1}\ \left( -1\right) ^{i}\ w_{i}\,x_{i}\
dx_{1}\wedge\cdots\wedge\widehat{dx_{i}}\wedge\cdots\wedge dx_{r+1},
\end{array}
\end{equation*}
and
\begin{equation*}
\begin{array}{l}
Gr_{\mathcal{F}^{\bullet}}^{p}(H^{r}(\mathbb{P}_{\mathbb{C}}^{r}\left(
\mathbf{w}\right) \mathbb{r}Z_{\infty},\mathbb{C}))\cong H^{r-p}(\mathbb{P}_{%
\mathbb{C}}^{r}\left( \mathbf{w}\right) ,\Omega_{\mathbb{P}_{\mathbb{C}%
}^{r}\left( \mathbf{w}\right) }^{p}(\log Z_{\infty})) \\
\  \\
\cong\dfrac{H^{0}(\mathbb{P}_{\mathbb{C}}^{r}\left( \mathbf{w}\right)
,\Omega_{\mathbb{P}_{\mathbb{C}}^{r}\left( \mathbf{w}\right)
}^{r}((r-p-1)Z_{\infty}))}{H^{0}(\mathbb{P}_{\mathbb{C}}^{r}\left( \mathbf{w}%
\right) ,\Omega_{\mathbb{P}_{\mathbb{C}}^{r}\left( \mathbf{w}\right)
}^{r}((r-p)Z_{\infty}))+\partial(H^{0}(\mathbb{P}_{\mathbb{C}}^{r}\left(
\mathbf{w}\right) ,\Omega_{\mathbb{P}_{\mathbb{C}}^{r}\left( \mathbf{w}%
\right) }^{r-1}((r-p)Z_{\infty})))}
\end{array}
\end{equation*}
($\partial$ denotes the corresponding differential operator). Since the map%
{\small
\begin{equation*}
\begin{array}{ccc}
H^{0}(\mathbb{P}_{\mathbb{C}}^{r}\left( \mathbf{w}\right) ,\Omega _{\mathbb{P%
}_{\mathbb{C}}^{r}\left( \mathbf{w}\right) }^{r}((r-p-1)Z_{\infty })) & \!%
\overset{\Phi}{\longrightarrow} & \!\mathbb{C}[x_{1},..,x_{r+1}]_{(r-p+1)d-%
\left( w_{1}+\cdots+w_{r+1}\right) }\smallskip \\
\  &  &  \\
\dfrac{g\cdot\Omega_{0}}{f^{r-p+1}}\  & \longmapsto & g
\end{array}
\end{equation*}
}defines an isomorphism, and
\begin{equation*}
H^{0}(\mathbb{P}_{\mathbb{C}}^{r}\left( \mathbf{w}\right) ,\Omega _{\mathbb{P%
}_{\mathbb{C}}^{r}\left( \mathbf{w}\right) }^{r}((r-p)Z_{\infty
}))+\partial(H^{0}(\mathbb{P}_{\mathbb{C}}^{r}\left( \mathbf{w}\right)
,\Omega_{\mathbb{P}_{\mathbb{C}}^{r}\left( \mathbf{w}\right)
}^{r-1}((r-p)Z_{\infty})))
\end{equation*}
has
\begin{equation*}
\left( \frac{\partial f}{\partial x_{1}},\ldots,\frac{\partial f}{\partial
x_{r+1}}\right) _{(r-p+1)d-\left( w_{1}+\ldots+w_{r+1}\right) }
\end{equation*}
as its image under $\Phi$ (see \cite[\S11]{BATYREV-COX}), we get
\begin{equation*}
Gr_{\mathcal{F}^{\bullet}}^{p}(H^{r}(\mathbb{P}_{\mathbb{C}}^{r}\left(
\mathbf{w}\right) \mathbb{r}Z_{\infty},\mathbb{C}))\cong M\left( f\right)
_{(r-p+1)d-\left( w_{1}+\ldots+w_{r+1}\right) }.
\end{equation*}
Using Hard Lefschetz Theorem one deduces the exact MHS-sequence:%
{\footnotesize
\begin{equation*}
0\!\rightarrow\! H^{r-2}(\mathbb{P}_{\mathbb{C}}^{r}\left( \mathbf{w}\right)
,\mathbb{C})\overset{\cong}{\longrightarrow}H^{r}(\mathbb{P}_{\mathbb{C}%
}^{r}\left( \mathbf{w}\right) ,\mathbb{C})\!\rightarrow H^{r}(\mathbb{P}_{%
\mathbb{C}}^{r}\left( \mathbf{w}\right) \mathbb{r}Z_{\infty},\mathbb{C}%
)\rightarrow H_{\text{prim}}^{r-1}(Z_{\infty },\mathbb{C})\rightarrow\!0,
\end{equation*}
}giving
\begin{equation*}
H_{\text{prim}}^{p,r-1-p}\left( Z_{\infty}\right) \cong M\left( f\right)
_{(r-p)d-\left( w_{1}+\ldots+w_{r+1}\right) }\cong M\left( f\right)
_{(p+1)d-\left( w_{1}+\ldots+w_{r+1}\right) }.
\end{equation*}
Formulae (\ref{HPQ-F}) follow from Lemma \ref{MIL-LEMMA} (iv), (v), and (\ref
{MIL-F-FDACH}).\hfill$\square$

\begin{lemma}
\label{HODGE-L}If $(X_{f},\mathbf{0})$ is an $r$-dimensional isolated
quasihomogeneous hypersurface singularity with $L$ as its link and $F$ as
its Milnor fiber$,$ then
\begin{equation*}
h^{p,q}\left( H^{r-1}\left( L,\mathbb{C}\right) \right) =0,\text{ whenever }%
p+q\neq r-1,\text{ }
\end{equation*}
and the ``non-trivial'' Hodge numbers of the cohomology groups of its link $L
$ are{\small
\begin{equation}
h^{p,r-1-p}\left( H^{r-1}\left( L,\mathbb{C}\right) \right)
=h^{r-p,p+1}\left( H^{r}\left( L,\mathbb{C}\right) \right)
=h^{p+1,r-p}\left( H^{r}\left( F,\mathbb{C}\right) \right)   \label{L-F}
\end{equation}
}for $p=0,1,\ldots ,r-1,$ and can be therefore read off from \emph{(\ref
{HPQ-F}).}
\end{lemma}

\noindent{}{}\textsc{Proof}. If $p+q\notin\{r-1,r+1\},$ then by \ref{COROL2}%
, \ref{COROL3} (iv) and \ref{MIL-LEMMA} (i), (iii), we obtain
\begin{equation*}
h^{r-p,r-q}\left( H^{r-1}\left( L,\mathbb{C}\right) \right) =h^{p,q}\left(
H^{r}\left( L,\mathbb{C}\right) \right) =h^{p,q}\left( H^{r-1}\left( L,%
\mathbb{C}\right) \right) =0,
\end{equation*}
because the corresponding Hodge numbers of $H^{r-1}\left( F,\mathbb{C}%
\right) $ and $H^{r}\left( F,\mathbb{C}\right) $ vanish, and $p+q<r-1$
(resp., $=r\,\ |>r+1$) iff $(r-p)+(r-q)>r+1$ (resp., $=r\,|<r-1$). On the
other hand, if $p+q\in\{r-1,r+1\},$ Cor. \ref{COROL2} gives:{\small
\begin{equation}
h^{p,q}(\!H^{r-1}(L,\mathbb{C})\!)-h^{p,q}(\!H^{r}(L,\mathbb{C}%
)\!)\!=\!h^{r-p,r-q}(\!H^{r}(F,\mathbb{C})\!)-h^{p,q}(\!H^{r}(F,\mathbb{C}%
)\!).  \label{TYPOS}
\end{equation}
}

{\small \noindent } \noindent {}If $p+q=r-1,$ the Hodge numbers $%
h^{p,q}(H^{r}(L,\mathbb{C}))=h^{r-p,r-q}(H^{r-1}(L,\mathbb{C}))$ vanish by
\ref{COROL3} (iv). Analogously, $h^{p,q}(H^{r-1}(L,\mathbb{C}))$ vanishes
whenever $p+q=r+1.$ Finally, (\ref{L-F}) follows from Lemma \ref{MIL-LEMMA}
(iii) and (\ref{TYPOS}).\hfill $\square $

\begin{proposition}
\label{BIG Prop}If $(X_{f},\mathbf{0})$ is an $r$-dimensional isolated
quasihomogeneous hypersurface singularity of degree $d$ w.r.t. the weights $%
w_{1},\ldots ,w_{r+1},$ and $L$ its link, then the E-polynomial $E\left(
X_{f}\mathbb{r}\{\mathbf{0}\};u,v\right) $ equals{\small
\begin{equation}
\left( uv-1\right) \ \left[ \sum_{p=0}^{r-1}\,\left( \left( uv\right)
^{p}+\left( -1\right) ^{r-1}\ h^{p,r-1-p}\left( H^{r-1}(L,\mathbb{C}\right)
)\ u^{p}\,v^{r-p-1}\right) \right]   \label{FORM-X}
\end{equation}
}and its coefficients are therefore computable in terms of $d$ and $%
w_{1},\ldots ,w_{r+1}$ via \emph{(\ref{L-F}) }and \emph{(\ref{HPQ-F}).}
\end{proposition}

\noindent {}\textsc{Proof}. Using (\ref{ISOL-X}) and Poincar\'{e} duality,
we obtain:
\begin{equation*}
h^{p,q}(H^{i}\left( L,\mathbb{C}\right) )=h^{p,q}(H^{i}\left( X_{f}\mathbb{r}%
\{\mathbf{0}\},\mathbb{C}\right) )=h^{d-p,d-q}(H_{c}^{2d-i}\left( X_{f}%
\mathbb{r}\{\mathbf{0}\},\mathbb{C}\right) ).
\end{equation*}
Hence,
\begin{equation}
E\left( X_{f}\mathbb{r}\{\mathbf{0}\};u,v\right) =\left( uv\right)
^{r}E\left( L;u^{-1},v^{-1}\right) .  \label{EXL}
\end{equation}
On the other hand, Proposition \ref{COROL3} gives{\footnotesize
\begin{equation*}
\begin{array}{l}
E\left( L;u,v\right) =\sum_{0\leq p,q\leq r}e^{p,q}(L)\ \,u^{p}\,v^{q}= \\
\, \\
=\sum\limits_{0\leq p,q\leq r}\ \left[ h^{p,q}(H^{0}(L))-h^{p,q}(H^{2r-1}(L))%
\right] \,\,u^{p}\,v^{q}+ \\
\, \\
+\left( -1\right) ^{r-1}\sum\limits_{0\leq p,q\leq r}\ \left[
\,h^{p,q}(H^{r-1}(L))-h^{p,q}(H^{r}(L))\right] \,\,u^{p}\,v^{q}= \\
\  \\
=\sum\limits_{0\leq p,q\leq r}\left[ h^{p,q}(H^{0}(L))-h^{p,q}(H^{2r-1}(L))%
\right] \,\ u^{p}\,v^{q}+ \\
\, \\
+\left( -1\right) ^{r-1}\sum\limits_{0\leq p,q\leq r}\left[
h^{p,q}(H^{r-1}(L))-h^{r-p,r-q}(H^{r-1}(L))\right] \,\ u^{p}\,v^{q}= \\
\  \\
=1-\left( uv\right) ^{r}+\left( -1\right) ^{r-1}\,[\sum\limits_{0\leq
p,q\leq r}\ h^{p,q}(H^{r-1}(L))\,]\,u^{p}\,v^{q}+ \\
\  \\
+\left( -1\right) ^{r}\,[\sum\limits_{0\leq p,q\leq r}\
h^{r-p,r-q}(H^{r-1}(L))\,\,u^{p}\,v^{q}]= \\
\, \\
=1-\left( uv\right) ^{r}+\left( -1\right) ^{r-1}\left[ \sum\limits
_{\substack{ 0\leq p,q\leq r-1  \\ 0\leq p+q\leq r-1}}\
h^{p,q}(H^{r-1}(L))\,\,u^{p}\,v^{q}\right] + \\
\  \\
+\left( -1\right) ^{r}\left[ \sum\limits_{\substack{ 1\leq p,q\leq r  \\ %
r+1\leq p+q\leq 2r-1}}\ h^{r-p,r-q}(H^{r-1}(L))\,\,u^{p}\,v^{q}\right] .
\end{array}
\end{equation*}
}(The terms containing coefficients $h^{p,q}(H^{r-1}(L))\,,$ with
$p+q=r,$ cancel out, as they occur in both summands). Since
$(X_{f},\mathbf{0})$ is an isolated quasihomogeneous hypersurface
singularity, we may use Lemma \ref {HODGE-L} to
write{\footnotesize
\begin{equation*}
\begin{array}{l}
E\left( L;u,v\right) =1-\left( uv\right) ^{r}+\left( -1\right) ^{r-1}\left[
\sum\limits_{\substack{ 0\leq p,q\leq r-1  \\ p+q=r-1}}\
h^{p,q}(H^{r-1}(L))\,\,u^{p}\,v^{q}\right] + \\
\  \\
+\left( -1\right) ^{r}\left[ \sum\limits_{\substack{ 1\leq p,q\leq r  \\ %
p+q=r+1}}\ h^{r-p,r-q}(H^{r-1}(L))\,\,u^{p}\,v^{q}\right] = \\
\  \\
=1-\left( uv\right) ^{r}+\left( -1\right) ^{r-1}\left[ \sum%
\limits_{p=0}^{r-1}\ h^{p,r-1-p}(H^{r-1}(L))\,\,u^{p}\,v^{r-1-p}\right] + \\
\  \\
+\left( -1\right) ^{r}\left[ \sum\limits_{p=0}^{r-1}\
h^{p,r-1-p}(H^{r-1}(L))\,\,u^{p+1}\,v^{r-p}\right] = \\
\  \\
=1-\left( uv\right) ^{r}+\left( -1\right) ^{r-1}\left[ \sum%
\limits_{p=0}^{r-1}\ h^{p,r-1-p}(H^{r-1}(L))\,\,u^{p}\,v^{r-1-p}\right]
\left( 1-uv\right) = \\
\  \\
=\left( 1-uv\right) \sum\limits_{p=0}^{r-1}\ (uv)^{p}+\left( -1\right)
^{r-1}\left( 1-uv\right) \,\left[ \sum\limits_{p=0}^{r-1}\
h^{p,r-1-p}(H^{r-1}(L))\,\,u^{p}\,v^{r-1-p}\right] = \\
\  \\
=\left( 1-uv\right) \left[ \sum\limits_{p=0}^{r-1}\ \left( (uv)^{p}+\left(
-1\right) ^{r-1}\,h^{p,r-1-p}(H^{r-1}(L))\,\,u^{p}\,v^{r-1-p}\right) \right]
\\
\
\end{array}
\end{equation*}
}Combining the last equality with (\ref{EXL})$,$ using
\begin{align*}
h^{p,r-1-p}(H^{r-1}(L,\mathbb{C}))& =\left( -1\right)
^{r-1}e^{p,r-1-p}\left( L\right) \\
& =\left( -1\right) ^{r-1}e^{r-1-p,p}\left( L\right) =h^{r-1-p,p}(H^{r-1}(L,%
\mathbb{C}))
\end{align*}
and substituting $r-p-1$ for $p,$ we deduce formula (\ref{FORM-X}).\hfill $%
\square $

\begin{remark}
\label{REM-RAT}(i) By (\ref{FORM-X}), $e(X_{f}\mathbb{r}\left\{ \mathbf{0}%
\right\} )=E\left( X_{f}\mathbb{r}\left\{ \mathbf{0}\right\} ;1,1\right)
=e(L)=0,$ which is also obvious from the fact, that $L$ is an \textit{odd}%
-dimensional differentiable manifold.\smallskip\ \newline
(ii) If the singularity $(X_{f},\mathbf{0})$ in \ref{BIG Prop} is, in
addition, a \textit{rational} singularity, then
\begin{equation*}
h^{0,r-1}(H^{r-1}(L))=h^{r-1,0}(H^{r-1}(L))=0.
\end{equation*}
(See the proof of Proposition 4.1 of \cite{DAIS-ROCZEN}.)\smallskip \newline
(iii) The defining polynomial (\ref{DEFF}) of an $\mathbf{A}_{n,\ell
}^{\left( r\right) }$-singularity is quasihomogeneous of degree $d=\text{lcm}%
\left( n+1,\ell \right) $ w.r.t. the weights $\left( \frac{d}{n+1},\frac{d}{%
\ell },\frac{d}{\ell },\ldots ,\frac{d}{\ell }\right) ,$ with Poincar\'{e}
polynomial
\begin{equation}
P_{M\left( f\right) }\left( t\right) =\left( 1+{\textstyle%
\sum\limits_{j=1}^{n-1}}t^{\frac{jd}{n+1}}\right) \,\left( 1+{\textstyle%
\sum\limits_{\kappa =1}^{\ell -2}}t^{\frac{\kappa d}{\ell }}\right) ^{r}
\label{POINC-SR}
\end{equation}
and Milnor number $\mu \left( f\right) =n\left( \ell -1\right) ^{r}$ (see (%
\ref{POIN-POL}) and (\ref{MILNOR-NUMBER})). Moreover, since $\mathbf{A}%
_{n,\ell }^{\left( r\right) }$'s are canonical, they are also rational
singularities.
\end{remark}

\noindent{}\textbf{Proof of Proposition \ref{PROP-FS}: }To produce formula (%
\ref{FIRST SUMMAND}) for the $E$-polynomial of $X_{n,\ell}^{(r)}\mathbb{r\{}%
\mathbf{0}\mathbb{\}}$, it suffices to evaluate (\ref{FORM-X}) via (\ref{L-F}%
), (\ref{HPQ-F}) and \ref{REM-RAT} (ii)-(iii) in terms of $n,\ell$ and $d.$
Since the function $\mathbf{a,}$ defined in (\ref{A-FUNCTION}), can be
expressed as
\begin{equation*}
\mathbf{a}\left( i\right) =\left\{
\begin{array}{ll}
1, & \text{if \thinspace}i\in\{0,\frac{d}{n+1},\frac{2d}{n+1},\ldots ,\frac{%
(n-1)d}{n+1}\} \\
0, & \text{otherwise,}
\end{array}
\right.
\end{equation*}
and since $\mathbf{b}\left( j\right) $, as defined in (\ref{B-FUNCTION}),
gives the coefficient of $t^{j}$ in the multinomial expansion of the second
factor of (\ref{POINC-SR}), we need the convolutional function (\ref
{CONVOLUTION}) in order to write the required dimensions as{\small
\begin{align*}
h^{p,r-1-p}\left( H^{r-1}\left( L,\mathbb{C}\right) \right) &
=h^{p+1,r-p}\left( H^{r}\left( F,\mathbb{C}\right) \right) \\
& =\dim_{\mathbb{C}}(M\left( f\right) _{d(p+1-\tfrac{1}{n+1}-\tfrac{r}{\ell}%
)})=\mathbf{c}(d(p+1-\tfrac{1}{n+1}-\tfrac{r}{\ell})),
\end{align*}
}and to end up to (\ref{FIRST SUMMAND}).\hfill$\square\bigskip$

\section{Desingularization and Theorem's Proof}

\noindent{}Next, using blow-ups of closed points, we shall construct
snc-resolutions for all $\mathbf{A}_{n,\ell}^{\left( r\right) }$%
-singularities for which either $\ell\,\left| \,n\right. $ or $\ell\,\left|
\,n\right. +1$. Let $X=X_{n,\ell}^{(r)}$ be their underlying spaces and
denote by
\begin{equation*}
Y_{\ell}^{\left( r-1\right) }:=\left\{ \left[ z_{1}:z_{2}:\ldots :z_{r+1}%
\right] \in\mathbb{P}_{\mathbb{C}}^{r}\ \left| \
\sum_{j=1}^{r+1}z_{j}^{\ell}=0\right. \right\}
\end{equation*}
the $(r-1)$-dimensional Fermat hypersurface of degree $\ell\geq2$ in the
projective space $\mathbb{P}_{\mathbb{C}}^{r}.$

\begin{proposition}
\label{DESINGULARIZATION}\emph{(i)} If $n+1\equiv 0\left( \emph{mod}\text{ }%
\ell \right) ,$ then there exists an snc- desingularization \ \ $\varphi :%
\widetilde{X}\longrightarrow X$ \ \ with discrepancy\textrm{\
\begin{equation}
K_{\widetilde{X}}-\varphi ^{\ast }\left( K_{X}\right) ={\textstyle%
\sum\limits_{i=1}^{m}}\ i\,\left( r-\ell \right) \ D_{i}
\label{DISCREPANCY1}
\end{equation}
where}
\begin{equation*}
D_{i}\cong \mathbb{P}(\mathcal{O}_{Y_{\ell }^{\left( r-2\right) }}\oplus
\mathcal{O}_{Y_{\ell }^{\left( r-2\right) }}\left( 1\right) ),\ \forall i,\
1\leq i\leq m-1,\text{ \ \emph{and} \ \ }D_{m}\cong Y_{\ell }^{\left(
r-1\right) }.
\end{equation*}
\emph{(ii)} If $n\equiv 0\left( \emph{mod}\text{ }\ell \right) ,$ then there
is an snc-desingularization $\varphi :\widetilde{X}\longrightarrow X$ \ with
discrepancy\textrm{\ {\small
\begin{equation}
K_{\widetilde{X}}-\varphi ^{\ast }\left( K_{X}\right) ={\textstyle%
\sum\limits_{i=1}^{m-1}}\ i\,\left( r-\ell \right) \ D_{i}+\left[ \left(
m-1\right) \,\ell \,\left( r-\ell \right) +\,\left( r-1\right) \right] \
D_{m}  \label{DISCREPANCY2}
\end{equation}
} where}
\begin{equation*}
D_{i}\cong \mathbb{P}(\mathcal{O}_{Y_{\ell }^{\left( r-2\right) }}\oplus
\mathcal{O}_{Y_{\ell }^{\left( r-2\right) }}\left( 1\right) ),\ \forall i,\
1\leq i\leq m-1,\text{ \ \emph{and} \ \ }D_{m}\cong \mathbb{P}_{\mathbb{C}%
}^{r-1}.
\end{equation*}
In both cases $D_{i}\cap D_{i+1}\cong Y_{\ell }^{\left( r-2\right) }$ for
all $i,1\leq i\leq m-1,$ and $D_{i}\cap D_{j}=\varnothing $ for all $%
i,j,1\leq i,j\leq m,$ with $\left| i-j\right| \neq 1.$ \emph{(}The number $m$
is defined as in \emph{(\ref{DEF-M})).}
\end{proposition}

\noindent\textsc{{}Proof}. Let $f$ be the polynomial (\ref{DEFF}) and $X$ $%
(=X_{f}=X_{n,\ell}^{(r)})$ the underlying space of the $\mathbf{A}_{n,\ell
}^{\left( r\right) }$-singularity.\bigskip

\noindent {}$\blacktriangleright $ \underline{\textbf{Construction of the
desingularization}}. Let $\pi :\mathbf{Bl}_{\mathbf{0}}(\mathbb{C}%
^{r+1})\longrightarrow \mathbb{C}^{r+1}$ be the blow up of $\mathbb{C}^{r+1}$
at the origin, with{\footnotesize
\begin{equation*}
\mathbf{Bl}_{\mathbf{0}}(\mathbb{C}^{r+1})=\left\{ \left( \left(
x_{1},..,x_{r+1}\right) ,\left[ t_{1}:\cdots :t_{r+1}\right] \right) \in
\mathbb{C}^{r+1}\times \mathbb{P}_{\mathbb{C}}^{r}\left|
\begin{array}{l}
\ x_{i}\,t_{j}=x_{j}\,t_{i}, \\
\ \forall i,j,\  \\
1\leq i,j\leq r+1
\end{array}
\right. \right\}
\end{equation*}
}$\mathcal{E}:=\pi ^{-1}\left( \mathbf{0}\right) =\left\{ \mathbf{0}\right\}
\times \mathbb{P}_{\mathbb{C}}^{r}$, and let $U_{i}$ denote the open set
given by $\left( t_{i}\neq 0\right) $. In terms of analytic coordinates,%
{\footnotesize
\begin{equation*}
U_{i}=\left\{ \!(\left( x_{1},..,x_{r+1}\right) ,(\xi _{1},..,\widehat{\xi
_{i}},..,\xi _{r+1}))\in \mathbb{C}^{r+1}\times \mathbb{C}^{r}\ \left|
\begin{array}{l}
\ x_{j}=x_{i}\,\xi _{j},\ \forall j, \\
j\in \left\{ 1,..,r+1\right\} \mathbb{r}\left\{ i\right\}
\end{array}
\!\right. \right\} \!,
\end{equation*}
}where $\xi _{j}=\frac{t_{j}}{t_{i}}$. Identifying $U_{i}$ with a copy of $%
\mathbb{C}^{r+1}$ w.r.t. the coordinates $x_{i}$, $\xi _{1}$, \ldots , $\widehat{%
\xi _{i}},\ldots ,\xi _{r+1}$, the restriction $\pi \left| _{U_{i}}\right. $
\ is given by mapping {\small
\begin{equation*}
\begin{array}{c}
\mathbb{C}^{r+1}\ni (x_{i},\xi _{1},..,\widehat{\xi _{i}},..,\xi _{r+1}) \\
\downarrow \ \cong \\
\begin{array}{c}
(\left( x_{i}\,\xi _{1},..,x_{i}\,\xi _{i-1},x_{i},x_{i}\,\xi
_{i+1},..,x_{i}\,\xi _{r+1}\right) ,[\xi _{1}:..:\underset{i\text{-th pos.}}{%
\underbrace{1}}:..:\xi _{r+1}])\in U_{i} \\
\downarrow \ \pi \left| _{U_{i}}\right. \\
\left( x_{i}\,\xi _{1},\ldots ,x_{i}\,\xi _{i-1},x_{i},x_{i}\,\xi
_{i+1},\ldots ,x_{i}\,\xi _{r+1}\right)
\end{array}
\end{array}
\end{equation*}
}Further, $\mathcal{E}_{i}:=\mathcal{E}\cap U_{i}$ is described as the
coordinate hyperplane $\left( x_{i}=0\right) $; i.e., the open cover $%
\left\{ U_{i}\right\} _{1\leq i\leq r+1}$ \ of $\mathbf{Bl}_{\mathbf{0}}(%
\mathbb{C}^{r+1})$ restricts to $\mathcal{E}$ to provide the standard open
cover of $\mathbb{P}_{\mathbb{C}}^{r}\ $by affine spaces $\mathbb{C}^{r+1}$,
with $\{\xi _{j}\}_{j\in \left\{ 1,\ldots ,r+1\right\} \mathbb{r}\left\{
i\right\} }$ being the analytic coordinates of $\mathcal{E}_{i}$.\smallskip

\noindent\textit{Notation. \ }To work with a more convenient notation we
define
\begin{equation*}
\mathbf{Bl}_{\mathbf{0}}(\mathbb{C}^{r+1})=\bigcup_{i=1}^{r+1}\ U_{i},\ \ \
\ \ U_{i}=\text{Spec}\left( \mathbb{C}\left[ y_{i,1},\ldots ,y_{i,r+1}\right]
\right) ,
\end{equation*}
by setting as coordinates for $U_{i}$'s:
\begin{equation*}
y_{i,k}:=\left\{
\begin{array}{ll}
x_{k}, & \text{for\ }i=k \\
\xi_{k}, & \text{for\ }i\neq k
\end{array}
\right.
\end{equation*}
\medskip

\noindent$\bullet$ \textsc{the first blow-up}\textbf{.} \ Blowing up $X$ at
the origin, we take the diagram\medskip\
\begin{equation*}
\begin{array}{ccccc}
\mathcal{E} & \subset & \mathbf{Bl}_{\mathbf{0}}(\mathbb{C}^{r+1}) &
\overset{\pi}{\longrightarrow} & \mathbb{C}^{r+1} \\
\cup &  & \cup &  & \cup \\
\mathcal{E}_{X}:=\mathcal{E}\cap\mathbf{Bl}_{\mathbf{0}}(X) & \subset &
\mathbf{Bl}_{\mathbf{0}}(X) & \overset{\pi\left| _{\text{restr.}}\right. }{%
\longrightarrow} & X
\end{array}
\end{equation*}
and consider the strict transform
\begin{equation*}
\mathbf{Bl}_{\mathbf{0}}(X)=\overline{\pi^{-1}(X\cap(\mathbb{C}^{r+1}\mathbb{%
r}\left\{ \mathbf{0}\right\} ))}=\overline{\pi^{-1}(X)\cap (\mathbf{Bl}_{%
\mathbf{0}}(\mathbb{C}^{r+1})\mathbb{r}\mathcal{E}))}
\end{equation*}
of $X$ in $\mathbb{C}^{r+1}$ under $\pi$, and the corresponding exceptional
divisor $\mathcal{E}_{X}$.\medskip

\noindent$\bullet$ \textsc{local description of}\textbf{\ }$\mathbf{Bl}_{%
\mathbf{0}}(X)$ \textsc{and} $\mathcal{E}_{X}.$ \ Pulling back $f$, we
get\medskip
\begin{equation*}
\pi^{\ast}(f)\left| _{U_{i}}\right. =x_{i}^{\ell}\ \widetilde{f_{i}}%
=y_{i,i}^{\ell}\ \widetilde{f_{i}}
\end{equation*}
with{\small $\ \widetilde{f_{i}}\left( y_{i,1},\ldots,y_{i,r+1}\right) =$\
\begin{equation*}
\left\{
\begin{array}{ll}
y_{1,1}^{(n+1)-\ell}+y_{1,2}^{\ell}+\cdots+y_{1,r+1}^{\ell}, & \text{%
if\thinspace\ }i=1 \\
\, &  \\
y_{i,1}^{n+1}\,y_{i,i}^{(n+1)-\ell}+y_{i,2}^{\ell}+\cdots+y_{i,i-1}^{\ell
}+1+y_{i,i+1}^{\ell}+\cdots+y_{i,r+1}^{\ell}, & \text{otherwise}
\end{array}
\right.
\end{equation*}
}

\noindent{}Locally,
\begin{equation*}
\mathbf{Bl}_{\mathbf{0}}(X)\left| _{U_{i}}\right. \cong\left\{ \left(
y_{i,1},\ldots,y_{i,r+1}\right) \in\mathbb{C}^{r+1}\ \left| \ \ \ \widetilde
{f_{i}}\left( y_{i,1},\ldots,y_{i,r+1}\right) =0\right. \right\} ,
\end{equation*}
and the equations for $\mathcal{E}_{X}\left| _{U_{i}}\right. $ read as
follows:{\small
\begin{align*}
\mathbf{Bl}_{\mathbf{0}}(X)\cap\mathcal{E}_{i} & =\mathcal{E}_{X}\left|
_{U_{i}}\right. \\
& \cong\left\{ \left( y_{i,1},\ldots,y_{i,r+1}\right) \in\mathbb{C}^{r+1}\
\left| \ \ y_{i,i}=\ \widetilde{f_{i}}\left( y_{i,1},\ldots
,y_{i,r+1}\right) =0\right. \right\} .
\end{align*}
}Thus, the only singular affine patch is $U_{1}=$ Spec$\left( \mathbb{C}%
\left[ y_{1,1},\ldots,y_{1,r+1}\right] \right) $ whenever $n>\ell.\medskip$

\noindent{}$\bullet$ \textsc{global description of}\textbf{\ }$\mathbf{Bl}_{%
\mathbf{0}}(X)$ \textsc{and} $\mathcal{E}_{X}.$ $\ $Passing to global
coordinates, we can write{\footnotesize
\begin{equation*}
\mathbf{Bl}_{\mathbf{0}}(X)=\left\{ \!\left( \left( x_{1},..,x_{r+1}\right) ,%
\left[ t_{1}:\cdots:t_{r+1}\right] \right) \in\mathbf{Bl}_{\mathbf{0}}(%
\mathbb{C}^{r+1})\left| x_{1}^{\left( n+1\right) -\ell }\,t_{1}^{\ell}+{%
\textstyle\sum\limits_{j=2}^{r+1}} t_{j}^{\ell}=0\right. \!\right\}
\end{equation*}
}and $\mathcal{E}_{X}$ equals{\footnotesize
\begin{equation*}
\left\{
\begin{array}{ll}
\left\{ \left( \mathbf{0},\left[ t_{1}:\cdots:t_{r+1}\right] \right) \in\
\left\{ \mathbf{0}\right\} \times\mathbb{P}_{\mathbb{C}}^{r}\ \left| \
\right. t_{1}^{\ell}+t_{2}^{\ell}+\cdots+t_{r+1}^{\ell}=0\right\} , & \text{%
if\thinspace\ }n=\ell-1 \\
\, &  \\
\left\{ \left( \mathbf{0},\left[ t_{1}:\cdots:t_{r+1}\right] \right) \in\
\left\{ \mathbf{0}\right\} \times\mathbb{P}_{\mathbb{C}}^{r}\ \left| \
\right. t_{2}^{\ell}+t_{3}^{\ell}+\cdots+t_{r+1}^{\ell}=0\right\} , & \text{%
otherwise}
\end{array}
\right. \medskip
\end{equation*}
}

{\footnotesize \noindent}$\bullet$ \textsc{the (fermat) singularity}\textbf{%
\ }$\mathbf{A}_{\ell-1,\ell}^{\left( r\right) }$ ($m=1$). Blowing up the
origin once, we achieve immediately the required desingularization, having
exceptional divisor $\mathcal{E}_{X}\cong Y_{\ell}^{\left( r-1\right)
}.\medskip$

\noindent$\bullet$ \textsc{the singularity} $\mathbf{A}_{\ell,\ell}^{\left(
r\right) }$ ($m=2$). In this case, $\mathbf{Bl}_{\mathbf{0}}(X)$ is smooth,
whereas $\mathcal{E}_{X}\subset\mathbf{Bl}_{\mathbf{0}}(X)$ has a singular,
ordinary $\ell$-fold point at
\begin{equation*}
Q=\left( \mathbf{0},[1:0:0:\cdots:0]\right) \in\mathcal{E}_{X}\left|
_{U_{1}}\right. .
\end{equation*}
To obtain an snc-resolution of the original singularity, we blow up once
more at $Q,$ and consider $\varphi=\pi_{1}\circ\pi_{2},$%
\begin{equation*}
\widetilde{X}=\mathbf{Bl}_{\,Q}(\mathbf{Bl}_{\mathbf{0}}(X))\overset{\pi_{2}%
}{\longrightarrow}\mathbf{Bl}_{\mathbf{0}}(X)\overset{\pi_{1}=\pi }{%
\longrightarrow}X\ .
\end{equation*}
\ The new exceptional divisor $D_{2}$ is a $\mathbb{P}_{\mathbb{C}}^{r-1}$,
and the strict transform $D_{1}$ of $\mathcal{E}_{X}$ is nothing but the (($%
r-1)$-dimensional) blow-up of $\mathcal{E}_{X}$ at $Q.$ Since $\mathcal{E}%
_{X}$ can be viewed as the projective cone $C^{\text{pr}}(Y_{\ell }^{\left(
r-2\right) })$ $\subset\mathbb{P}_{\mathbb{C}}^{r}$ over the Fermat
hypersurface $Y_{\ell}^{\left( r-2\right) }\subset\mathbb{P}_{\mathbb{C}%
}^{r-1}$ with $[1:0:\cdots:0]$ as its vertex, blowing up $[1:0:\cdots:0],$
the diagram{\small
\begin{equation*}
\begin{array}{ccccc}
\mathbb{P(\mathcal{O}}_{\mathbb{P}_{\mathbb{C}}^{r-1}}\oplus \mathbb{%
\mathcal{O}}_{\mathbb{P}_{\mathbb{C}}^{r-1}}(1)\mathbb{)} & \cong & \mathbf{%
Bl}_{[1:0:\cdots:0]}(\mathbb{P}_{\mathbb{C}}^{r}) & \longrightarrow &
\mathbb{P}_{\mathbb{C}}^{r} \\
&  & \cup &  & \cup \\
&  & \mathbf{Bl}_{[1:0:\cdots:0]}(C^{\text{pr}}(Y_{\ell}^{\left( r-2\right)
})) & \longrightarrow & C^{\text{pr}}(Y_{\ell}^{\left( r-2\right) })\cong%
\mathcal{E}_{X}
\end{array}
\end{equation*}
}yields the isomorphism
\begin{equation*}
\mathbf{Bl}_{[1:0:\cdots:0]}(C^{\text{pr}}(Y_{\ell}^{\left( r-2\right)
}))\cong\mathbb{P(\mathcal{O}}_{Y_{\ell}^{\left( r-2\right) }}\oplus\mathbb{%
\mathcal{O}}_{Y_{\ell}^{\left( r-2\right) }}(1)\mathbb{)}.
\end{equation*}
Hence, $D_{1}$ is a $\mathbb{P}_{\mathbb{C}}^{1}$-bundle of rank $2$ over $%
Y_{\ell}^{\left( r-2\right) }$ meeting $D_{2}$ along
\begin{equation*}
\left( D_{1}\cdot D_{2}\right) \left| _{D_{1}}\right. =\mathbb{P(\mathcal{O}}%
_{Y_{\ell}^{\left( r-2\right) }}(1)\mathbb{)}\cong Y_{\ell}^{\left(
r-2\right) }
\end{equation*}
(see Fig. 1).

\begin{figure}[h]
\begin{center}
\includegraphics[width=271.1875pt,height=247.125pt]{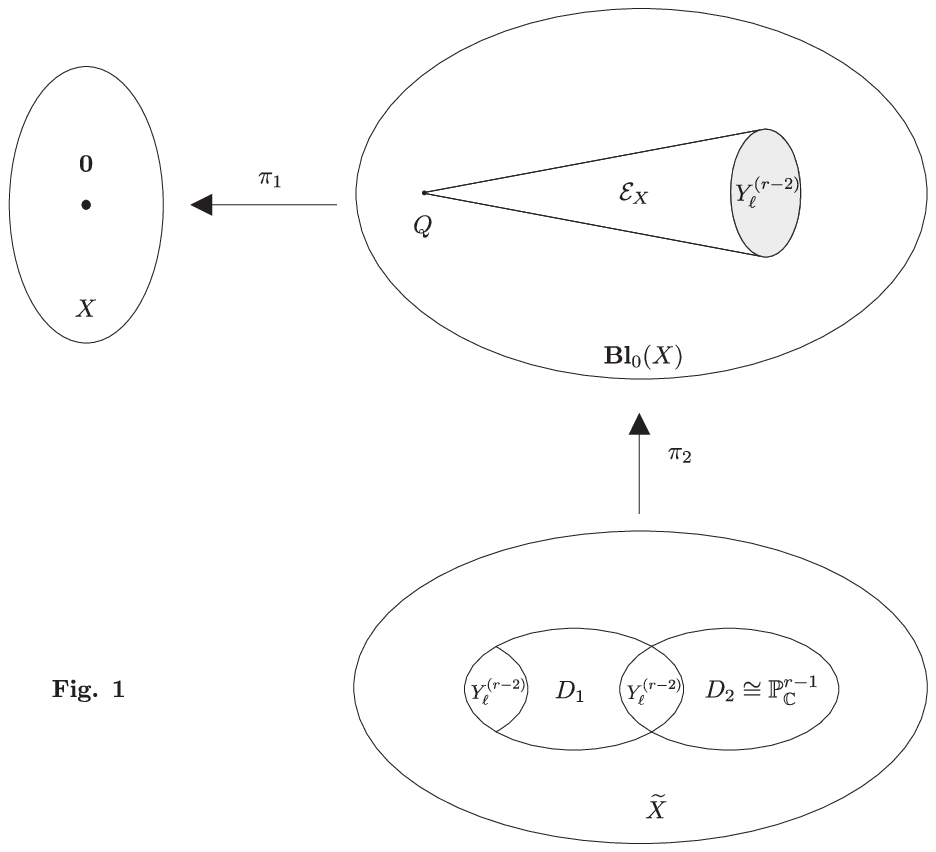}
\end{center}
\end{figure}

\noindent {}$\bullet $ \textsc{singularities} $\mathbf{A}_{n,\ell }^{\left(
r\right) }$ $\ $\textsc{with}\textbf{\ }$n>\ell ,$ $m\geq 2,$ \textsc{and
either} $\ell \,\left| \,n\right. $\textsc{\ or} $\ell \,\left| \,n\right.
+1 $. Locally, these singularities can be reduced successively to one of the
above types as follows:\medskip
\begin{equation*}
\fbox{$
\begin{array}{l}
\mathbf{A}_{n,\ell }^{\left( r\right) }\rightsquigarrow \mathbf{A}_{n-\ell
,\ell }^{\left( r\right) }\rightsquigarrow \mathbf{A}_{n-2\ell ,\ell
}^{\left( r\right) }\rightsquigarrow \cdots \rightsquigarrow \mathbf{A}%
_{2\ell -1,\ell }^{\left( r\right) }\rightsquigarrow \mathbf{A}_{\ell
-1,\ell }^{\left( r\right) }\rightsquigarrow \mathbf{A}_{-1,\ell }^{\left(
r\right) } \\
\  \\
\ (\text{if\ \ }n+1\equiv 0\left( \text{mod}\text{ }\ell \right) ) \\
\  \\
\begin{array}{l}
\mathbf{A}_{n,\ell }^{\left( r\right) }\rightsquigarrow \mathbf{A}_{n-\ell
,\ell }^{\left( r\right) }\rightsquigarrow \mathbf{A}_{n-2\ell ,\ell
}^{\left( r\right) }\rightsquigarrow \cdots \rightsquigarrow \mathbf{A}%
_{2\ell ,\ell }^{\left( r\right) }\rightsquigarrow \mathbf{A}_{\ell ,\ell
}^{\left( r\right) }\rightsquigarrow \mathbf{A}_{0,\ell }^{\left( r\right)
}\rightsquigarrow \mathbf{A}_{0,\ell }^{\left( r\right) } \\
\  \\
(\text{if\ \ }n\equiv 0\left( \text{mod}\text{ }\ell \right) )
\end{array}
\end{array}
$}\medskip
\end{equation*}
(Each ``$\rightsquigarrow "$ denotes the result of a local blow-up, and $%
\mathbf{A}_{-1,\ell }^{\left( r\right) }$, $\mathbf{A}_{0,\ell }^{\left(
r\right) }$ stand for ``smooth charts''). But also globally, $\varphi :%
\widetilde{X}\longrightarrow X$ \ is decomposed into just $m$ blow-ups
\begin{equation}
\begin{array}{l}
\widetilde{X}=X_{m}\overset{\pi _{m}}{\longrightarrow }X_{m-1}\overset{\pi
_{m-1}}{\longrightarrow }\cdots \overset{\pi _{3}}{\longrightarrow }X_{2}%
\overset{\pi _{2}}{\longrightarrow }X_{1}\overset{\pi _{1}}{\longrightarrow }%
X_{0}=X \\
\, \\
X_{i}:=\mathbf{Bl}_{\,Q_{i}}(\mathbf{Bl}_{\,Q_{i-1}}(\,\cdots \,(\mathbf{Bl}%
_{\,Q_{1}}(X)))),\text{ \ }\forall i,\ 1\leq \text{ }i\leq m,
\end{array}
\label{COMPOSITION}
\end{equation}
of $m$ ``separated'' points $Q_{1}=\mathbf{0},\ Q_{2}=\left( \mathbf{0}%
,[1:0:0:\cdots :0]\right) ,\ldots ,Q_{m}$, in the sense, that all the
appearing exceptional divisors are prime (by construction) and, in addition,
if $E_{1}=\mathcal{E}_{X},E_{2},\ldots ,E_{m}$ are the exceptional loci of $%
\pi _{1},\pi _{2},\ldots ,\pi _{m},$ respectively, the singular point $Q_{i}$
is resolved by $\pi _{i}$ and the (possibly existing) new singular point $%
Q_{i+1}$ is \textit{not} contained in the strict transforms of $%
E_{1},E_{2},\ldots ,E_{i-1}$ under $\pi _{i}$. Thus, defining the divisor $%
D_{i}$ to be the strict transform of $E_{i}$ under $\pi _{i+1}\circ \pi
_{i+2}\circ \cdots \circ \pi _{m-1}\circ \pi _{m}$ on $\widetilde{X}$, we
obtain the intersection graph of Figure 2 with $D_{i}\cap D_{i+1}\cong
Y_{\ell }^{\left( r-2\right) }$.

\begin{figure}[h]
\begin{center}
\includegraphics[width=11.1cm,height=86.5pt]{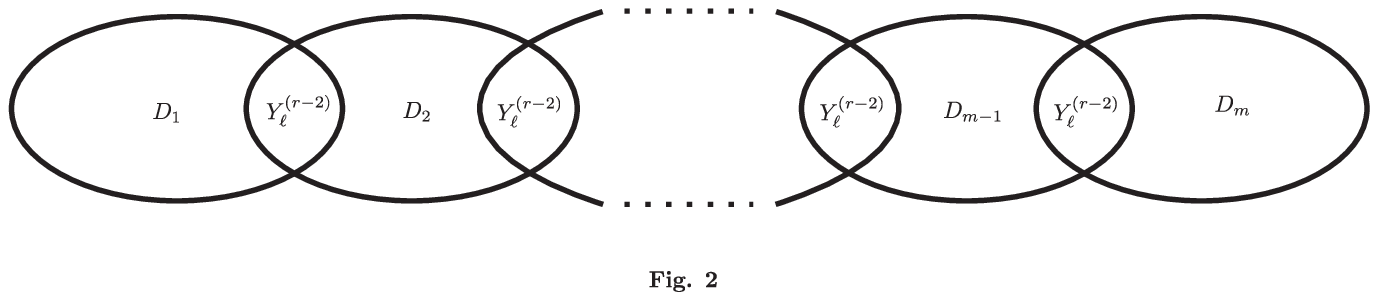}
\end{center}
\end{figure}

\noindent {}$\blacktriangleright $ \underline{\textbf{Computation of the
discrepancy coefficients}}. Consider the Poincar\'{e} residue map
\begin{equation*}
\text{Res}_{X}:H^{0}(\mathbb{C}^{r+1},\omega _{\mathbb{C}^{r+1}}(X))%
\longrightarrow H^{0}(X,\omega _{X}),
\end{equation*}
where $\omega _{X}=\mathcal{O}_{X}(K_{X})=(\Omega _{X}^{r})^{\vee \,\vee
}\subset \Omega _{\mathbb{C}\left( X\right) /\mathbb{C}}^{r}$. The rational
canonical differential
\begin{equation*}
\mathfrak{s:}=\text{Res}_{X}\left( \frac{dx_{1}\wedge dx_{2}\wedge
dx_{3}\wedge \cdots \wedge dx_{r+1}}{f}\right) =\frac{dx_{2}\wedge
dx_{3}\wedge \cdots \wedge dx_{r+1}}{\left( \partial f\,/\,\partial
x_{1}\right) }
\end{equation*}
can be viewed as a (local) generator of $H^{0}(X,\omega _{X})$. Assume that $%
n\neq \ell $ and that you have performed the first blow-up of $X$ at $%
\mathbf{0.}$ Then the new singularity (if any) on $\mathbf{Bl}_{\mathbf{0}%
}(X)$ will belong to $\mathcal{E}_{X}\left| _{U_{1}}\right. $. For this
reason, to find the discrepancy coefficient w.r.t. $\pi :\mathbf{Bl}_{%
\mathbf{0}}(X)\longrightarrow X,$ it suffices to compare $\mathfrak{s}$ with
the rational canonical differential
\begin{equation*}
\overline{\mathfrak{s}}:=\frac{dy_{1,2}\wedge dy_{1,3}\wedge \cdots \wedge
dy_{1,r+1}}{(\partial \widetilde{f_{1}}\,/\,\partial y_{1,1})}\in \Omega _{%
\mathbb{C}\left( U_{1}\right) /\mathbb{C}}^{r}\ .
\end{equation*}
($U_{1}$ is non-singular with local coordinates $y_{1,2},\ldots ,y_{1,r+1}$
at any point $P$ for which $\partial \widetilde{f_{1}}(P)\,/\,\partial
y_{1,1}\neq 0$). In $U_{1}$ we have $x_{1}=y_{1,1}$ and
\begin{equation*}
x_{j}=x_{1}\,\xi _{j}=y_{1,1}\,y_{1,j},\text{ \ for all \ }j\in
\{2,3,...,r+1\}.
\end{equation*}
Hence,
\begin{equation}
\begin{array}{l}
dx_{2}\wedge dx_{3}\wedge \cdots \wedge dx_{r+1} \\
\  \\
=y_{1,1}^{r-1}\ (y_{1,1}\left( dy_{1,2}\wedge dy_{1,3}\wedge \cdots \wedge
dy_{1,r+1}\right) + \\
\  \\
+\sum\limits_{i=2}^{r+1}\ \left( -1\right) ^{i}\ y_{1,i}\ dy_{1,2}\wedge
\cdots \wedge \widehat{dy_{1,i}}\wedge \cdots \wedge dy_{1,r+1})
\end{array}
\label{DF1}
\end{equation}
and
\begin{equation}
\partial f\,/\,\partial x_{1}=\left( n+1\right) \,x_{1}^{n}=\left(
n+1\right) \,\,y_{1,1}^{n}=\left( \tfrac{n+1}{n+1-\ell }\right) \
y_{1,1}^{\ell }\ (\partial \widetilde{f_{1}}\,/\,\partial y_{1,1})
\label{DIF2}
\end{equation}
On the other hand,
\begin{equation*}
d\widetilde{f_{1}}=\left( n+1-\ell \right) \ y_{1,1}^{n-\ell
}\,dy_{1,1}+\ell \,\left( y_{1,2}^{\ell -1}\,dy_{1,2}+\cdots
+y_{1,r+1}^{\ell -1}\,dy_{1,r+1}\right) =0
\end{equation*}
if and only if
\begin{equation}
dy_{1,1}=-\frac{\ell }{n+1-\ell }\ y_{1,1}^{\ell -n}\ \left( y_{1,2}^{\ell
-1}\,dy_{1,2}+\cdots +y_{1,r+1}^{\ell -1}\,dy_{1,r+1}\right)  \label{DIF3}
\end{equation}
Substituting the expression (\ref{DIF3}) for $dy_{1,1}$ into the right-hand
side of (\ref{DF1}), we obtain
\begin{equation}
\begin{array}{l}
dx_{2}\wedge dx_{3}\wedge \cdots \wedge dx_{r+1} \\
\  \\
=\left( -\tfrac{\ell }{n+1-\ell }\ y_{1,1}^{r-1+\ell -n}\ (y_{1,2}^{\ell
}+\cdots +y_{1,r+1}^{\ell })+y_{1,1}^{r}\right) \ dy_{1,2}\wedge \cdots
\wedge dy_{1,r+1}
\end{array}
\label{DIF4}
\end{equation}
Combining now (\ref{DIF4}) with $y_{1,2}^{\ell }+\cdots +y_{1,r+1}^{\ell
}=-y_{1,1}^{(n+1)-\ell }$ and (\ref{DIF2}), we get
\begin{equation}
\mathfrak{s}=\frac{y_{1,1}^{r}\ dy_{1,2}\wedge dy_{1,3}\wedge \cdots \wedge
dy_{1,r+1}}{y_{1,1}^{\ell }\ (\partial \widetilde{f_{1}}\,/\,\partial
y_{1,1})}=y_{1,1}^{r-\ell }\ \overline{\mathfrak{s}}  \label{DIF5}
\end{equation}
The equality (\ref{DIF5}) shows that the discrepancy coefficient of $%
\mathcal{E}_{X}$ with respect to $\pi :\mathbf{Bl}_{\mathbf{0}%
}(X)\longrightarrow X$ equals $r-\ell $. Using the notation introduced in (%
\ref{COMPOSITION}), one proves analogously that
\begin{equation}
K_{X_{i}}-\pi _{i}^{\ast }\left( K_{X_{i-1}}\right) =\left( r-\ell \right)
\,E_{i},\ \ \forall i,\ 1\leq \text{ }i\leq m-1.  \label{K1}
\end{equation}
Moreover,
\begin{equation}
K_{X_{m}}-\pi _{m}^{\ast }\left( K_{X_{m-1}}\right) =\left\{
\begin{array}{ll}
\left( r-\ell \right) D_{m}, & \text{if \ }\ell \left| n+1\right. \\
\, & \  \\
\left( r-1\right) D_{m}, & \text{if \ }\ell \left| n\right.
\end{array}
\right.  \label{K2}
\end{equation}
Note that if $\ell \left| n\right. ,$ then we have to pass through $\mathbf{A%
}_{\ell ,\ell }^{\left( r\right) }.$ The additional blow-up which resolves
the singularity of the exceptional locus (so that $\varphi :\widetilde{X}%
\longrightarrow X$ fulfills the snc-condition) has a \textit{smooth} point
on the $r$-\textit{fold} as its centre. Consequently, the discrepancy
coefficient of $D_{m}=D_{\frac{n}{\ell }+1}$ equals $r-1$ (see \cite[p. 187]
{GR-H})$.$ Now (\ref{K1}) gives:
\begin{equation}
\begin{array}{l}
K_{\widetilde{X}}-\varphi ^{\ast }\left( K_{X}\right) = \\
\, \\
{\textstyle}\sum\limits_{i=1}^{m-1}\left( \pi _{i+1}\circ \pi _{i+2}\circ
\cdots \circ \pi _{m}\right) ^{\ast }\left( \left( r-\ell \right)
E_{i}\right) +\left[ K_{X_{m}}-\pi _{m}^{\ast }\left( K_{X_{m-1}}\right)
\right]
\end{array}
\label{k3}
\end{equation}
Since{\small
\begin{equation}
\left( \pi _{i+1}\circ \pi _{i+2}\circ \cdots \circ \pi _{m}\right) ^{\ast
}\left( \,E_{i}\right) =\left\{
\begin{array}{ll}
{\textstyle}\sum\limits_{j=i}^{m}D_{j}, & \text{if \ }\ell \left| n+1\right.
\\
\, & \, \\
{\textstyle}\sum\limits_{j=i}^{m-1}D_{j}+\ell \,D_{m}, & \text{if \ }\ell
\left| n\right.
\end{array}
\right.  \label{K4}
\end{equation}
}for all $i,$ $1\leq i\leq m-1,$ the formulae (\ref{DISCREPANCY1}) and (\ref
{DISCREPANCY2}) follow from (\ref{K2}), (\ref{k3}) and (\ref{K4}).\hfill $%
\square $

\begin{remark}
\noindent {}\label{FINAL REMARK}(i) If $n+1\equiv 0\left( \text{mod}\text{ }%
\ell \right) $ and $r=\ell ,$ then $\varphi :\widetilde{X}\longrightarrow X$
is crepant.\smallskip \newline
(ii) Obviously,
\begin{equation}
E(\mathbb{P}_{\mathbb{C}}^{r-1};u,v)={\textstyle\sum\limits_{p=0}^{r-1}}%
\left( uv\right) ^{p}.\newline
\label{E-PRSP}
\end{equation}
(iii) To complete the catalogue of the $E$-polynomials of our exceptional
divisors, it suffices to find out those of $Y_{\ell }^{\left( r-2\right) }$
(or, equivalently, of $Y_{\ell }^{\left( r-1\right) }),$ as we have
\begin{equation}
E(\mathbb{P}(\mathcal{O}_{Y_{\ell }^{\left( r-2\right) }}\oplus \mathcal{O}%
_{Y_{\ell }^{\left( r-2\right) }}\left( 1\right) );u,v)=E(Y_{\ell }^{\left(
r-2\right) };u,v)\cdot \left( 1+uv\right)   \label{E-BUNDLE}
\end{equation}
(iv) According to the classical Lefschetz Hyperplane Theorem, the Fermat
hypersurface $Y_{\ell }^{\left( r-1\right) }$ has ``non-trivial'' Hodge $%
\left( p,q\right) $-numbers only if $p+q=r-1.$ Next lemma expresses them by
means of the non-central Eulerian numbers of generalized factorials (as
defined in \S 1 \textsf{(d)}), and can be easily proven, e.g., by
determining the $\chi _{y}$-characteristic of $Y_{\ell }^{\left( r-1\right) }
$ via Riemann-Roch Theorem (see \cite[\S2]{HIRZEBRUCH1}), or, alternatively,
by writing down the exact sequences involving the cohomology groups of $%
\mathbb{P}_{\mathbb{C}}^{r}$ and $Y_{\ell }^{\left( r-1\right) }$ with
coefficients taken from the twisted sheaves $\Omega _{\mathbb{P}_{\mathbb{C}%
}^{r}}^{p}(-\ell )$ and $\Omega _{Y_{\ell }^{\left( r-1\right) }}^{p}(-\ell )
$, respectively. (Note that both proofs are valid for \textit{any} smooth
hypersurface of degree $\ell .$ On the other hand, the formula for the Euler
number is simpler and can be derived directly by evaluating the highest
Chern class of $Y_{\ell }^{\left( r-1\right) }$ and applying Gauss-Bonnet
Theorem; see, e.g., \cite[p. 152]{DIMCA2}.)
\end{remark}

\begin{lemma}
The Hodge numbers of the $(r-1)$-dimensional Fermat hypersurface $Y_{\ell
}^{\left( r-1\right) }$ of degree $\ell \geq 2$ are given by the formula
\begin{equation*}
h^{p,q}(Y_{\ell }^{\left( r-1\right) })=\left\{
\begin{array}{ll}
\mathfrak{S}\left( r,p+1\ \left| \ \ell -1,p\right. \right) +\delta
_{2p,r-1}, & \text{\emph{if} \ \ }p+q=r-1 \\
\  & \, \\
\delta _{p,q}, & \text{\emph{if} \ \ }p+q\neq r-1
\end{array}
\right.
\end{equation*}
Hence, $Y_{\ell }^{\left( r-1\right) }$ has $E$-polynomial{\small
\begin{align}
E(Y_{\ell }^{\left( r-1\right) };u,v)& ={\textstyle\sum\limits_{0\leq
p,q\leq r-1}}\,\left( -1\right) ^{p+q}\,h^{p,q}(Y_{\ell }^{\left( r-1\right)
})\,u^{p}v^{q}  \label{E-FERMAT} \\
& ={\textstyle\sum\limits_{p=0}^{r-1}}\,u^{p}\,\left[ v^{p}+\left( -1\right)
^{r-1}\mathfrak{S}\left( r,p+1\ \left| \ \ell -1,p\right. \right) \,v^{r-1-p}%
\right]   \notag
\end{align}
}and Euler number
\begin{equation}
\begin{array}{l}
e(Y_{\ell }^{\left( r-1\right) })=\left[ {\textstyle\sum\limits_{p=0}^{r-1}}%
\,\left( -1\right) ^{r-1}\mathfrak{S}\left( r,p+1\ \left| \ \ell -1,p\right.
\right) \right] +r= \\
\  \\
=\tfrac{1}{\ell }\left( \left( 1-\ell \right) ^{r+1}-1\right) +r+1
\end{array}
\label{e-FERMAT}
\end{equation}
\end{lemma}

\noindent {}\textbf{Proof of Theorem \ref{MAIN}.} (i) If $n+1\equiv 0\left(
\text{mod}\text{ }\ell \right) ,$ then Proposition \ref{DESINGULARIZATION}
and (\ref{E-STR}) give:{\small
\begin{equation*}
\begin{array}{l}
E_{\text{str}}(X;u,v)-E(X\mathbb{r}\{\mathbf{0}\};u,v)= \\
\, \\
{\textstyle}{\textstyle\sum\limits_{i=1}^{m}}\tfrac{\left( uv-1\right)
E(D_{i}^{\circ };u,v)}{\left( uv\right) ^{i(r-\ell )+1}-1}+\left[ {\textstyle%
}\sum\limits_{i=1}^{m-1}\tfrac{\left( uv-1\right) ^{2}\
E(D_{\{i,i+1\}}^{\circ };u,v)}{\left( \left( uv\right) ^{i(r-\ell
)+1}-1\right) \left( \left( uv\right) ^{(i+1)(r-\ell )+1}-1\right) }\right]
\\
\
\end{array}
\end{equation*}
}But
\begin{equation*}
E(D_{1}^{\circ };u,v)=\left\{
\begin{array}{ll}
E(Y_{\ell }^{\left( r-1\right) };u,v) & \text{if }\ell =n+1\text{ (i.e., }m=1%
\text{)} \\
E(Y_{\ell }^{\left( r-2\right) };u,v)\cdot uv, & \text{otherwise,}
\end{array}
\right.
\end{equation*}
\begin{equation*}
\,E(D_{i}^{\circ };u,v)=E(Y_{\ell }^{\left( r-2\right) };u,v)\cdot \left(
uv-1\right) ,\,\forall i,\,i\in \{2,\ldots ,m-1\},\,
\end{equation*}
(by (\ref{E-BUNDLE})), and (for $m\geq 2$) $E(D_{m}^{\circ };u,v)=E(Y_{\ell
}^{\left( r-1\right) };u,v)-E(Y_{\ell }^{\left( r-2\right) };u,v),$
\begin{equation*}
E(D_{\{i,i+1\}}^{\circ };u,v)=E(Y_{\ell }^{\left( r-2\right)
};u,v),\,\forall i,\,i\in \{1,\ldots ,m-1\}.
\end{equation*}
Consequently, for $m\geq 2,$ the difference $E_{\text{str}}(X;u,v)-E(X%
\mathbb{r}\{\mathbf{0}\};u,v)$ equals{\small
\begin{equation*}
\left( uv-1\right) \,E(Y_{\ell }^{\left( r-2\right) };u,v)\,\left[ \tfrac{uv%
}{\left( uv\right) ^{r-\ell +1}-1}+{\textstyle}\sum_{i=2}^{m-1}\tfrac{uv-1}{%
\left( uv\right) ^{i(r-\ell )+1}-1}-\tfrac{uv-1}{\left( uv\right) ^{m(r-\ell
)+1}-1}\right]
\end{equation*}
}
\begin{equation*}
+\tfrac{\left( uv-1\right) E(Y_{\ell }^{\left( r-1\right) };u,v)}{\left(
uv\right) ^{m(r-\ell )+1}-1}+{\textstyle}\sum\limits_{i=1}^{m-1}\tfrac{%
\left( uv-1\right) ^{2}\ E(Y_{\ell }^{\left( r-2\right) };u,v)}{\left(
\left( uv\right) ^{i(r-\ell )+1}-1\right) \left( \left( uv\right)
^{(i+1)(r-\ell )+1}-1\right) },
\end{equation*}
leading to the desired formula via (\ref{E-FERMAT}). [For $m=1,$ the
computation is staightforward]. Passing to the limit of $E_{\text{str}%
}(X;u,v)$, for $u,v\rightarrow 1,$ and taking (\ref{e-STR}) and (\ref
{e-FERMAT}) into account, one obtains the corresponding formula for the
string-theoretic Euler number $e_{\text{str}}(X).\medskip $

\smallskip \noindent {}\noindent {(ii)}\thinspace ~ If $n\equiv 0\left(
\text{mod}\text{ }\ell \right) ,$ then \ref{DESINGULARIZATION} (ii) and (\ref
{E-STR}) give analogously
\begin{equation*}
\begin{array}{l}
E_{\text{str}}(X;u,v)-E(X\mathbb{r}\{\mathbf{0}\};u,v)= \\
\, \\
{\textstyle}{\textstyle\sum\limits_{i=1}^{m-1}}\tfrac{\left( uv-1\right)
E(D_{i}^{\circ };u,v)}{\left( uv\right) ^{i(r-\ell )+1}-1}+\tfrac{uv-1}{%
\left( uv\right) ^{\left( m-1\right) \,\ell \,\left( r-\ell \right) +\,r}-1}%
\,E(D_{m}^{\circ };u,v) \\
\, \\
+\left[ {\textstyle}\sum\limits_{i=1}^{m-2}\tfrac{\left( uv-1\right) ^{2}\
E(D_{\{i,i+1\}}^{\circ };u,v)}{\left( \left( uv\right) ^{i(r-\ell
)+1}-1\right) \left( \left( uv\right) ^{(i+1)(r-\ell )+1}-1\right) }\right]
\\
\  \\
+\tfrac{\left( uv-1\right) ^{2}\ E(D_{\{m-1,m\}}^{\circ };u,v)}{\left(
\left( uv\right) ^{(m-1)(r-\ell )+1}-1\right) \left( \left( uv\right)
^{\left( m-1\right) \,\ell \,\left( r-\ell \right) +\,r}-1\right) }
\end{array}
\end{equation*}
Since $D_{1}^{\circ },D_{2}^{\circ },\ldots ,D_{m-1}^{\circ }$ are as in
(i), and
\begin{equation*}
E(D_{m}^{\circ };u,v)=E(\mathbb{P}_{\mathbb{C}}^{r-1};u,v)-E(Y_{\ell
}^{\left( r-2\right) };u,v),
\end{equation*}
\begin{equation*}
E(D_{\{i,i+1\}}^{\circ };u,v)=E(Y_{\ell }^{\left( r-2\right) };u,v),\text{
for all \ }i\in \{1,\ldots ,m-1\},\,
\end{equation*}
we obtain by (\ref{E-PRSP}) (\ref{E-BUNDLE}):{\small
\begin{equation*}
\begin{array}{l}
E_{\text{str}}(X;u,v)-E(X\mathbb{r}\{\mathbf{0}\};u,v)= \\
\, \\
\begin{array}{l}
\left( uv-1\right) \,E(Y_{\ell }^{\left( r-2\right) };u,v)\,\left[ \tfrac{uv%
}{\left( uv\right) ^{r-\ell +1}-1}+{\textstyle}\sum\limits_{i=2}^{m-1}\tfrac{%
uv-1}{\left( uv\right) ^{i(r-\ell )+1}-1}-\tfrac{uv-1}{\left( uv\right)
^{\left( m-1\right) \,\ell \,\left( r-\ell \right) +\,r}-1}\right] \\
\, \\
+\tfrac{\left( uv-1\right) \left( {\textstyle}\sum_{p=0}^{r-1}\left(
uv\right) ^{p}\right) }{\left( uv\right) ^{\left( m-1\right) \,\ell \,\left(
r-\ell \right) +\,r}-1}+\left[ {\textstyle}\sum\limits_{i=1}^{m-2}\tfrac{%
\left( uv-1\right) ^{2}\ E(Y_{\ell }^{\left( r-2\right) };u,v)}{\left(
\left( uv\right) ^{i(r-\ell )+1}-1\right) \left( \left( uv\right)
^{(i+1)(r-\ell )+1}-1\right) }\right] \\
\  \\
+\tfrac{\left( uv-1\right) ^{2}\ E(Y_{\ell }^{\left( r-2\right) };u,v)}{%
\left( \left( uv\right) ^{(m-1)(r-\ell )+1}-1\right) \left( \left( uv\right)
^{\left( m-1\right) \,\ell \,\left( r-\ell \right) +\,r}-1\right) }
\end{array}
\\
\,
\end{array}
\end{equation*}
}The string-theoretic Euler number is examined as in (i).\hfill$\square $

\section{Some global geometric examples}

{}\noindent {}The $E_{\text{str}}$-function of a complex $r$-fold $V$ with
only $k$ $\ $isolated log-terminal singularities $Q_{1},Q_{2},..,Q_{k}$
equals:
\begin{equation}
E_{\text{str}}\left( V;u,v\right) =E\left( V;u,v\right) \ +\ \sum_{i=1}^{k}\
\left( E_{\text{str}}\left( \left( V,Q_{i}\right) ;u,v\right) -1\right)
\label{ESTR-GLOBAL}
\end{equation}
In particular, a simple closed formula for the string-theoretic Euler number
$e_{\text{str}}$ can be easily built whenever $V$ is a (global) complete
intersection in a projective space, equipped with prescribed singularities
belonging to the class under consideration.

\begin{proposition}
Let $V=V_{\left( d_{1},d_{2},\ldots ,d_{N-r}\right) }$ $\ $be an $r$%
-dimensional complete intersection of multidegree $\left( d_{1},d_{2},\ldots
,d_{N-r}\right) $ in $\mathbb{P}_{\mathbb{C}}^{N}$ having only $k$ isolated
singularities $Q_{1},\ldots ,Q_{k}$ of types $\mathbf{A}_{n_{1},\ell
_{1}}^{\left( r\right) },$ $\ldots ,\mathbf{A}_{n_{k},\ell _{k}}^{\left(
r\right) }$ with either $\ell _{i}\left| n\right. _{i}$ or $\ell _{i}\left|
n\right. _{i}+1,$ for all $i=1,...,k.$ Then its string-theoretic Euler
number equals{\footnotesize
\begin{equation*}
e_{\text{\emph{str}}}(V)\!=\!\left[ \!\tbinom{N+1}{r}+{\textstyle%
\sum\limits_{\nu =1}^{r}}\,\left( {\textstyle\sum\limits_{1\leq j_{1}\leq
\cdots \leq j_{\nu }\leq N-r}}d_{j_{1}}\,\cdots \,d_{j_{\nu }}\right) \!%
\right] \left( {\textstyle\prod\limits_{j=1}^{N-r}}\ d_{j}\right) \ \!+
\end{equation*}
}

{\footnotesize
\begin{equation}
+\ \ \sum\limits_{i=1}^{k}\ \left[ e_{\text{\emph{str}}}\left(
V,Q_{i}\right) +\left( -1\right) ^{r+1}\,n_{i}\left( \ell _{i}-1\right)
^{r}-1\right]   \label{ESTR-CI}
\end{equation}
}where $e_{\text{\emph{str}}}\left( Y,Q_{i}\right) ,\ i=1,...,k,$ are
computable via Theorem\emph{\ \ref{MAIN}.}
\end{proposition}

\noindent\textsc{Proof}. By a small deformation of $V$ one can always obtain
a non-singular complete intersection $V^{\prime}$ in $\mathbb{P}_{\mathbb{C}%
}^{N}$ having multidegree $\left( d_{1},d_{2},\ldots,d_{N-r}\right) $. Using
a standard technique which involves the Mayer-Vietoris sequence (cf.
\cite[Ch. 5, Cor. 4.4 (ii) ]{DIMCA2}) one shows easily that
\begin{equation*}
e\left( V\right) =e\left( V^{\prime}\right) \ +\ \left( -1\right) ^{r+1}\
\sum\limits_{i=1}^{k}\text{\ }\left[ \text{Milnor number of }\left(
V,Q_{i}\right) \right] \ .
\end{equation*}
The Euler number of $V^{\prime}$ can be computed again by evaluating the
highest Chern class of $V^{\prime}$ at its fundamental cycle (cf. Chen-Ogiue
\cite[Thm. 2.1]{CH-O}), and is expressible by the closed formula:{\small
\begin{equation*}
\begin{array}{c}
e\left( V^{\prime}\right) =\left[ \!\binom{N+1}{r}+\sum\limits_{\nu=1}^{r}%
\,(-1)^{\nu}\,\binom{N+1}{r-\nu}\,\left( \sum\limits_{1\leq
j_{1}\leq\cdots\leq j_{\nu}\leq N-r}d_{j_{1}}\,\cdots\,d_{j_{\nu}}\right) \!%
\right] \left( \prod\limits_{j=1}^{N-r}\ d_{j}\right) .
\end{array}
\end{equation*}
}(\ref{ESTR-CI}) follows clearly from (\ref{ESTR-GLOBAL}).\hfill$\square$

\begin{examples}
Let us now apply (\ref{ESTR-CI}) for some well-known hypersurfaces and
complete intersections.\medskip \newline
\textbf{(i)} Generalizing Hirzebruch's method of constructing a singular
quintic with $126$ nodes (\cite[p. 762]{HIRZEBRUCH3}), Werner defines in
\cite[pp. 216-217]{WERNER} a hypersurface $V\subset $ $\mathbb{P}_{\mathbb{C}%
}^{4}$ of degree $5$ by homogenizing a three-dimensional affine complex
variety of the form
\begin{equation*}
\left\{ \left( z_{1},z_{2},z_{3},z_{4}\right) \in \mathbb{C}^{4}\ \left| \
\right. \lambda \,g_{1}(z_{1},z_{2})-g_{2}(z_{3},z_{4})=0\right\} ,\ \lambda
\in \mathbb{C}^{\ast },
\end{equation*}
where $\{g_{i}=0\},$ $i=1,2,$ are plane quintic curves having the three axes
and a circumscribed conic (about the corresponding coordinate triangle) as
their irreducible components (see Fig. 3). Since each of these curves has $3$
$\mathbf{D}_{4}$-singularities, $V$ (after homogenization) will have $3^{2}=9
$ singularities of type $\mathbf{A}_{2,3}^{\left( 3\right) }.$ This means
that
\begin{equation*}
e_{\text{str}}(V)=-200+9\cdot \left( 9+2^{4}-1\right) =16
\end{equation*}
In fact, $e_{\text{str}}(V)=e(\widetilde{V})=16,$ where $\widetilde{V}%
\rightarrow V$ is the crepant desingularization of $V$ arising from a single
simultaneous blow-up of the $9$ singularities (cf. \ref{FINAL REMARK} (i)). $%
\widetilde{V}$ is obviously a $3$-dimensional Calabi-Yau manifold.
\begin{figure}[h]
\begin{center}
\includegraphics[width=127.9375pt,height=91.625pt]{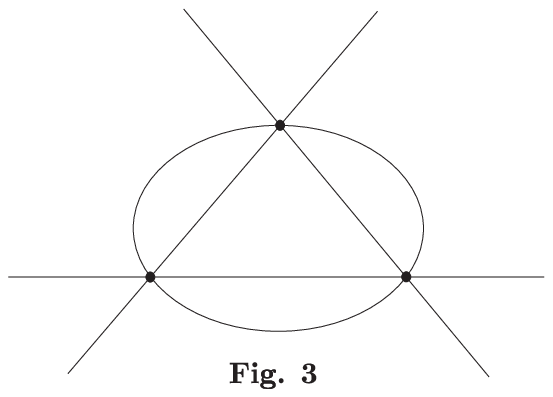}
\end{center}
\end{figure}
\newline
\bigskip \noindent {}\textbf{(ii)} The $\left( N-1\right) $-dimensional
Goryunov's quartics \cite{GORYUNOV}:{\small
\begin{equation*}
V_{\kappa }:=\left\{ \left[ z_{1}:..:z_{N+1}\right] \in \mathbb{P}_{\mathbb{C%
}}^{N}\left| 2\left( \kappa +1\right) {\textstyle\sum\limits_{1\leq i<j\leq
N+1}}z_{i}^{2}z_{j}^{2}+\kappa \left( {\textstyle\sum\limits_{1\leq j\leq
N+1}}z_{j}^{2}\right) ^{2}=\ 0\right. \right\}
\end{equation*}
}($N\geq \kappa ,$ $N\geq 3,$ $\kappa \geq 0$) have $2^{\kappa }\binom{N+1}{%
\kappa +1}$ $\mathbf{A}_{1}$-singularities ($\mathbf{A}_{1,2}^{\left(
N-1\right) }$-singularities, in our notation), and string-theoretic Euler
number
\begin{align*}
e_{\text{str}}(V_{\kappa })& =\tfrac{1}{4}\left( \left( -3\right)
^{N+1}-1\right) +N+1+ \\
& +2^{\kappa }\tbinom{N+1}{\kappa +1}\left[ \left( \tfrac{1}{N-2}\left(
\tfrac{1}{2}\left( \left( -1\right) ^{N}-1\right) +N\right) +\left(
-1\right) ^{N}-1\right) \right] .
\end{align*}
Note that, e.g., for $N=5,$ the string-theoretic index of the underlying
space of each of the singularities is $3>1,$ whereas the string-theoretic
index ind$_{\text{str}}(V_{\kappa })$ of $V_{\kappa }$ can be equal to $1$
(for $\kappa \in \left\{ 0,1,3,4\right\} $).\bigskip\ \newline
\textbf{(iii)} The $\left( n-2\right) $-dimensional Segre-Kn\"{o}rrer
complete intersection of two quadrics
\begin{equation*}
V:=\left\{ \mathbf{z}=\left[ z_{1}:z_{2}:\ldots :z_{n+1}\right] \in \mathbb{P%
}_{\mathbb{C}}^{n}\ \left| \ \ ^{t}\mathbf{z}\ M\ \mathbf{z}=\ ^{t}\mathbf{z}%
\ M^{\prime }\ \mathbf{z\ }=\ 0\right. \right\} ,\ \ n\geq 4,
\end{equation*}
where $M$ and $M^{\prime }$ are the $\left( n+1\right) \times \left(
n+1\right) $-matrices:
\begin{equation*}
M=\left(
\begin{array}{cccccc}
0 & 0 & \cdots  & \cdots  & 0 & 1 \\
0 & 0 & \cdots  & \cdots  & 1 & 0 \\
\vdots  & \vdots  & \cdots  & \cdots  & \vdots  & \vdots  \\
0 & 0 & \cdots  & \cdots  & 0 & 0 \\
0 & 1 & \cdots  & \cdots  & 0 & 0 \\
1 & 0 & \cdots  & \cdots  & 0 & 0
\end{array}
\right) \ ,\ M^{\prime }=\left(
\begin{array}{cccccc}
0 & 0 & 0 & \cdots  & 0 & 0 \\
0 & 0 & 0 & \cdots  & 0 & 1 \\
\vdots  & \vdots  & \cdots  & \cdots  & 1 & \vdots  \\
0 & 0 & 0 & \cdots  & 0 & 0 \\
0 & 0 & 1 & \cdots  & 0 & 0 \\
0 & 1 & 0 & \cdots  & 0 & 0
\end{array}
\right) ,
\end{equation*}
has $Q=[1:0:\cdots :0:0]$ as single isolated point which is of type $\mathbf{%
A}_{n}$ (i.e., $\mathbf{A}_{n,2}^{\left( n-2\right) }$ in our notation, see
\cite[p. 48]{KNOERRER}). According to (\ref{ESTR-CI}), the string-theoretic
Euler number of $V$ equals
\begin{align*}
e_{\text{str}}(V)& =\sum\limits_{\nu =0}^{n-2}\,(-1)^{\nu }\,2^{\nu +2}\,%
\tbinom{n+1}{\nu +3}\,\left( \nu +1\right) +e_{\text{str}}\left( V,Q\right)
+\left( -1\right) ^{n-1}\,n-1 \\
& =n-1+e_{\text{str}}\left( V,Q\right) ,
\end{align*}
with
\begin{equation*}
e_{\text{str}}\left( V,Q\right) =
\begin{cases}
\frac{(n-1)^{2}}{n^{2}-3n-2}, & \text{if }n\text{ odd} \\[6pt]
\,\frac{(n-2)(n+1)}{n(n-4)+(n-2)}, & \text{if }n\text{ even}
\end{cases}
\end{equation*}
For $n\leq 15,$ $e_{\text{str}}(V)$ takes the following values:
\begin{equation*}
\begin{tabular}{|c|c|c|c|c|c|c|c|c|c|c|c|c|}
\hline
$n$ & $4$ & $5$ & $6$ & $7$ & $8$ & $9$ & $10$ & $11$ & $12$ & $13$ & $14$ &
$15$ \\ \hline
$e_{\text{str}}(V)$ & $
\begin{array}{c}
\  \\
8 \\
\
\end{array}
$ & $6$ & $\frac{27}{4}$ & $\frac{96}{4}$ & $\frac{160}{19}$ & $\frac{120}{13%
}$ & $\frac{175}{17}$ & $\frac{480}{43}$ & $\frac{648}{53}$ & $\frac{105}{8}$
& $\frac{539}{38}$ & $\frac{1344}{89}$ \\ \hline
\end{tabular}
\end{equation*}
\newline
\textbf{(iv)}\thinspace ~Werner's $3$-dimensional complete intersection of a
cubic and two quadrics
\begin{equation*}
V:=\left\{ \left[ z_{1}:z_{2}:...:z_{7}\right] \in \mathbb{P}_{\mathbb{C}%
}^{6}\ \left| \begin{aligned} {\textstyle\sum\limits_{i=1}^{4}}
z_{i}^{3}&={\textstyle\sum\limits_{j=2}^{7}} z_{j}^{2} \\
&={\textstyle\sum\limits_{i=1}^{3}}
i\,z_{i+1}^{2}+{\textstyle\sum\limits_{j=4}^{6}} \left( j-3\right)
\,z_{j+1}^{2}=\ 0 \end{aligned}\right. \right\}
\end{equation*}
has $4$ singularities of type $\mathbf{A}_{2,3}^{\left( 3\right) }$ at the
points $\left[ 0:0:0:0:\pm 1:\pm \sqrt{-2}:1\right] $ and $18$~singularities
of type $\mathbf{A}_{1,2}^{\left( 3\right) }$ (i.e., nodes) at the points
\begin{align*}
& \bigl[-\zeta _{3}^{j}:1:0:0:\pm \sqrt{-1}:0:0\big],\ \ \big[-\zeta
_{3}^{j}:0:1:0:0:\pm \sqrt{-1}:0\big], \\
& \bigl[-\zeta _{3}^{j}:0:0:1:0:0:\pm \sqrt{-1}\big],
\end{align*}
$j=1,2,3,$ where $\zeta _{3}$ is a primitive third root of unity (see
\cite[pp. 221--222]{WERNER}). Its string-theoretic Euler number equals
\begin{equation*}
e_{\text{str}}(V)=-144+4\cdot \left( 9+2^{4}-1\right) +18\cdot 2=-12=e(%
\widetilde{V}),
\end{equation*}
where $\widetilde{V}$ is a Calabi--Yau threefold which arises after a
crepant desingularization of $V$ coming from the simultaneous (usual)
blow-up of the $9$ $\mathbf{A}_{2,3}^{\left( 3\right) }$-singularities and
an appropriate \textit{small, projective} resolution of the $18$
nodes.\medskip \newline
\textbf{(v)}\thinspace ~ Let $V=V_{1}\cap V_{2}\cap \cdots \cap
V_{N-r}\subset \mathbb{P}_{\mathbb{C}}^{N}$ be a complete intersection of
Fermat hypersurfaces
\begin{equation*}
V_{i}=\left\{ \left[ z_{1}:\ldots :z_{N+1}\right] \in \mathbb{P}_{\mathbb{C}%
}^{N}\ \left| \ \ {\textstyle\sum\limits_{j=1}^{N+1}}b_{ij}\,z_{j}^{d}=\
0\right. \right\} ,\quad 1\leq i\leq N-r,
\end{equation*}
of degree $d,$ $2\leq d\leq r,$ and assume that $V$ is $r$-dimensional,
i.e.,
\begin{equation*}
\text{rank}\left( (b_{ij})_{1\leq i\leq N-r,1\leq j\leq N+1}\right) =N-r.
\end{equation*}
Further, consider the map
\begin{equation*}
\Phi _{d}:\mathbb{P}_{\mathbb{C}}^{N}\longrightarrow \mathbb{P}_{\mathbb{C}%
}^{N},\,\,\,\left[ z_{1}:\ldots :z_{N+1}\right] \longmapsto \left[
z_{1}^{d}:\ldots :z_{N+1}^{d}\right] =\left[ \xi _{1}:\ldots :\xi _{N+1}%
\right] .
\end{equation*}
$\Phi _{d}$ displays $\mathbb{P}_{\mathbb{C}}^{N}$ as a $d^{N}$-sheeted
ramified covering of itself, branched along the coordinate axes $\{\xi
_{j}=0\}$. On the other hand,
\begin{equation*}
\Phi _{d}\left( V_{i}\right) =\left\{ \left[ \xi _{1}:..:\xi _{N+1}\right]
\in \mathbb{P}_{\mathbb{C}}^{N}\ \left| \ \ {\textstyle\sum%
\limits_{j=1}^{N+1}}b_{ij}\,\xi _{j}=\ 0\right. \right\} ,\quad 1\leq i\leq
N-r,
\end{equation*}
and $\Phi _{d}\left( V\right) \cong \mathbb{P}_{\mathbb{C}}^{r}\subset
\mathbb{P}_{\mathbb{C}}^{N}$. Now if
\begin{equation*}
\mathcal{L}_{j}:=\{\xi _{j}=0\}\cap \Phi _{d}\left( V\right) \subset \mathbb{%
P}_{\mathbb{C}}^{r},\quad 1\leq j\leq N+1,
\end{equation*}
denote by $\mathcal{M}\left( \mathbb{P}_{\mathbb{C}}^{N}\right) =\mathbb{C}%
\left( z_{2}/z_{1},\ldots ,z_{N+1}/z_{1}\right) $ the rational function
field of $\mathbb{P}_{\mathbb{C}}^{N},$ and let
\begin{equation*}
\mathcal{M}\left( \mathbb{P}_{\mathbb{C}}^{N}\right) \left( \sqrt[d]{\frac{%
\psi _{2}}{\psi _{1}}},\,\ldots ,\,\sqrt[d]{\frac{\psi _{N+1}}{\psi _{1}}}%
\right)
\end{equation*}
be the Kummer extension of $\mathcal{M}\left( \mathbb{P}_{\mathbb{C}%
}^{N}\right) $ determined by adjoining ``$d$-th roots of ratios'', where $%
\psi _{j}$ is the linear form defining the hyperplane $\mathcal{L}_{j}$.
This is an abelian extension with Galois group $\left( \mathbb{Z\,}/\,d%
\mathbb{Z}\right) ^{N}$. The variety $V$ can be thought of as the
normalization of $\mathbb{P}_{\mathbb{C}}^{N}$ w.r.t. this field, as being
the total space of the $d^{N}$-sheeted covering
\begin{equation*}
\Phi _{d}\left| _{V}\right. :V\longrightarrow \mathbb{P}_{\mathbb{C}}^{r}
\end{equation*}
of $\mathbb{P}_{\mathbb{C}}^{r},$ branched along the $\mathcal{L}_{j}$'s.
The hyperplane arrangement
\begin{equation*}
\mathfrak{L}:=\bigcup_{j=1}^{N+1}\,\mathcal{L}_{j}=\left\{ {\prod_{j=1}^{N+1}%
}\,\psi _{j}=0\right\} \subset \mathbb{P}_{\mathbb{C}}^{r}
\end{equation*}
admits a natural stratification {\small
\begin{equation*}
\begin{array}{c}
\mathfrak{L}=\mathfrak{L}^{\left( 1\right) }\!\supset \!\mathfrak{L}^{\left(
2\right) }=\!\!\!{\textstyle\bigcup\limits_{1\leq j_{1}<j_{2}\leq N+1}}%
\mathcal{L}_{j_{1},j_{2}}\supset \cdots \supset \mathfrak{L}^{\left(
r\right) }=\!\!\!{\textstyle\bigcup\limits_{1\leq j_{1}<j_{2}<\cdots
<j_{r}\leq N+1}}\mathcal{L}_{j_{1},j_{2},\ldots ,j_{r}}
\end{array}
\end{equation*}
} where
\begin{equation*}
\mathcal{L}_{j_{1},j_{2},\ldots ,j_{k}}:=\mathcal{L}_{j_{1}}\cap \mathcal{L}%
_{j_{2}}\cap \cdots \cap \mathcal{L}_{j_{k}}\cong \mathbb{P}_{\mathbb{C}%
}^{r-k}\subset \mathbb{P}_{\mathbb{C}}^{r},\ \ \ 1\leq k\leq r.
\end{equation*}
($\mathfrak{L}^{\left( r\right) }$ consists of the \textit{points} of $%
\mathfrak{L},$ $\mathfrak{L}^{\left( r-1\right) }$ consists of the \textit{%
lines} of $\mathfrak{L},$ etc). Let us now define
\begin{equation*}
t_{i}:=t_{i}\left( 0\right) :=\#\left\{ \begin{aligned}{}& \text{elements of
}\mathfrak{L}^{\left( r\right) }\text{ (i.e., points of }\mathfrak{
L}\text{) } \\ &\text{contained in exactly }i\text{ hyperplanes of
}\mathfrak{L} \end{aligned}\right\}
\end{equation*}
and, in general,
\begin{equation*}
t_{i}\left( \kappa \right) :=\#\left\{ \begin{aligned}{}& \text{elements of
}\mathfrak{L}^{\left( r-\kappa\right) }\text{ contained } \\ &\text{in
exactly }i\text{ hyperplanes of }\mathfrak{L} \end{aligned}\right\} ,\quad
0\leq \kappa \leq r.
\end{equation*}
$\mathfrak{L}$ is called a \textit{point arrangement} if
\begin{equation*}
t_{i}\left( \kappa \right) =0,\text{ \ for all \ }i>r-\kappa \text{ \ and
for all \ }\kappa \in \{1,\ldots ,r-2\}.
\end{equation*}
The $V$'s defined by means of point arrangements have at most isolated
singularities; more precisely, by analogy with the two-dimensional case (cf.
\cite{HIRZEBRUCH2}), $V$ inherits exactly $d^{N-i}$ isolated singularities
over each point of $\mathfrak{L}$ contained in $i\geq r+1$ hyperplanes. In
particular, for point arrangements $\mathfrak{L}$ within $\mathbb{P}_{%
\mathbb{C}}^{r}$, for which
\begin{equation*}
t_{i}=0,\quad \forall \,i,\quad i\geq r+2,
\end{equation*}
all singularities of $V$ have to be $\mathbf{A}_{d-1,d}^{\left( r\right) }$
-singularities. In this case, formula (\ref{ESTR-CI}) reads as follows:
\begin{align}
e_{\text{str}}(V)={}& \left[ \sum\limits_{\nu =0}^{r}\,(-1)^{\nu }\,\tbinom{%
N+1}{r-\nu }\,\tbinom{N-r+\nu -1}{\nu }\,d^{\nu +N-r}\right] +  \notag \\
& {}+\ t_{r+1}\cdot d^{N-r-1}\cdot \biggl(\tfrac{1}{r-d+1}\left[ \tfrac{1}{d}%
(\left( 1-d\right) ^{r+1}-1)+r+1\right]   \notag \\
& {}+\left( -1\right) ^{r+1}\,\left( d-1\right) ^{r+1}-1\biggl).
\label{HYPERFORM}
\end{align}
For $r=3$, several combinatorial properties of hyperplane arrangements in $%
\mathbb{P}_{\mathbb{C}}^{3},$ as well as properties of birational geometry
of the resulting coverings, have been studied by Hunt \cite{HUNT}. As far as
point arrangements are concerned (with $t_{4}\geq 1,$ $t_{5}=t_{6}=0$) there
are some interesting and aesthetically pleasing examples, given by the facet
planes of certain regular (platonic) and semiregular (archimedean) solids
(see Fig. 4). For these point arrangements, formula (\ref{HYPERFORM})
gives:\medskip
\begin{equation*}
\begin{tabular}{|c|c|c|c|c|c|}
\hline
Solids & $N$ & $
\begin{array}{c}
\  \\
t_{3} \\
\
\end{array}
$ & $t_{4}=\frac{1}{4}[\tbinom{N+1}{3}-t_{3}]$ & $
\begin{array}{c}
\  \\
e_{\text{str}}(V) \\
\  \\
(\text{for}\ d=2)
\end{array}
$ & $
\begin{array}{c}
\  \\
e_{\text{str}}(V) \\
\  \\
(\text{for}\ d=3)
\end{array}
$ \\ \hline\hline
\textbf{A} & $
\begin{array}{c}
5\
\end{array}
$ & $8$ & $3$ & $12$ & $\allowbreak \allowbreak 72$ \\ \hline
\textbf{B} & $7$ & $8$ & $12$ & $64$ & $\allowbreak -324$ \\ \hline
\textbf{C} & $7$ & $32$ & $6$ & $-32$ & $-4212$ \\ \hline
\textbf{D}, \textbf{E} & $13$ & $256$ & $27$ & $-111\,616$ & $\allowbreak
-68\,496\,840$ \\ \hline
\textbf{F} & $13$ & $208$ & $39$ & $-99\,328$ & $-62\,828\,136$ \\ \hline
\end{tabular}
\end{equation*}

\begin{figure}[h]
\begin{center}
\includegraphics[width=275.625pt,height=225.25pt]{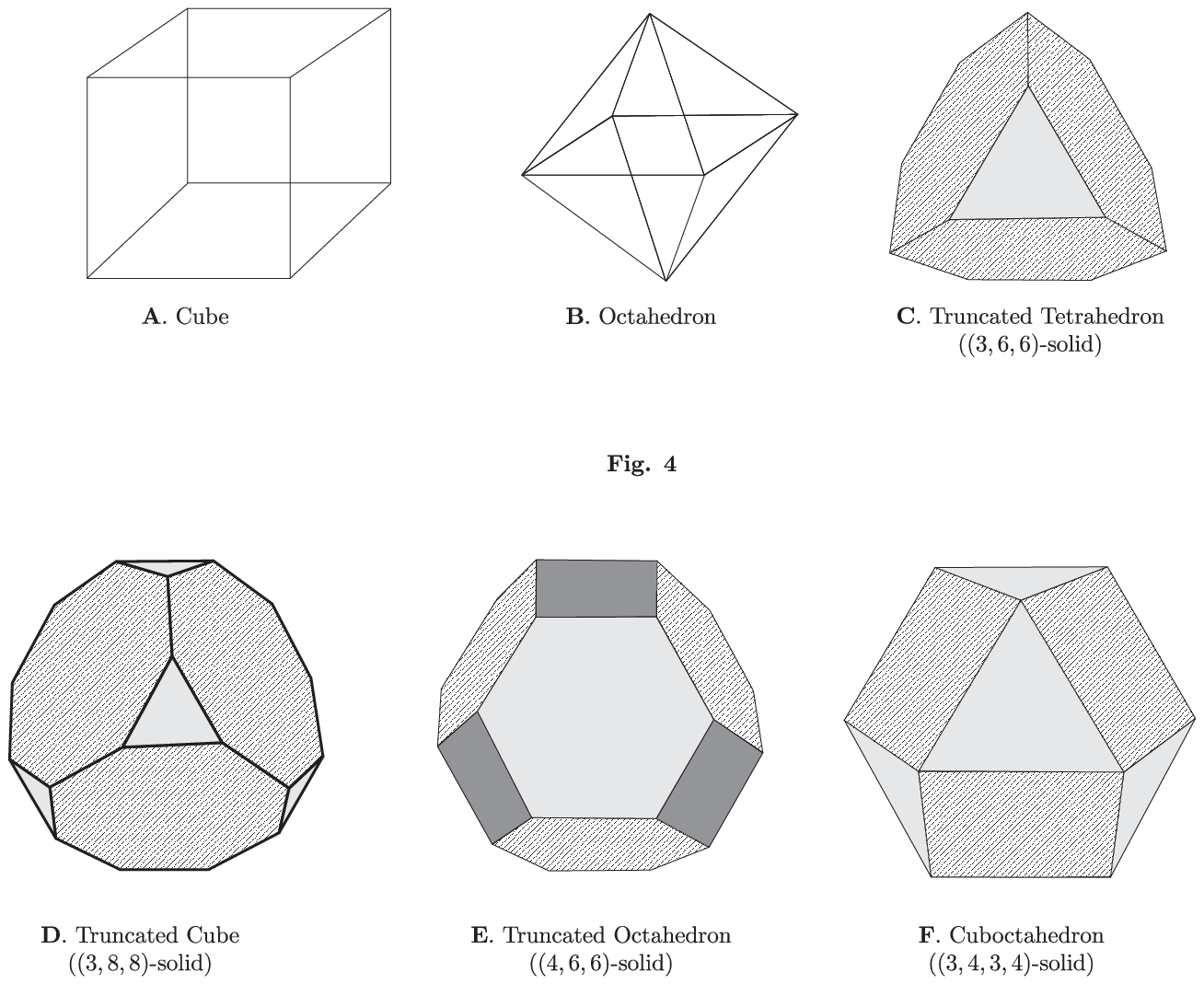}
\end{center}
\end{figure}

\noindent {}Examples \textbf{A} (with $d=3)$ and \textbf{B} (with $d=2$)
were first mentioned by Hirzebruch \cite[pp. 764--765]{HIRZEBRUCH3}, who
used them to construct $3$-dimensional Calabi-Yau manifolds $\widetilde{V}$
with Euler number $72$ (resp., $64$) by a ``big'' (resp. ``small'',
projective) crepant resolution of the $9$ (resp., $96$) singularities of $V$
(cf. the remarks in \cite[p.~219]{WERNER}).\medskip
\end{examples}

\noindent \textbf{Aknowledgements.} The author would like to
express his gratitude to Nobuyuki Kakimi (University of Tokyo) who
informed him about some ``missed'' extra factors of the
discrepancy coefficients (cf. (\ref {DISCREPANCY1}),
(\ref{DISCREPANCY2})) in a previous version of the paper. The
particular form of these coefficients led to the counterexamples
which are mentioned in \ref{CONJREM}.

\end{document}